\documentclass{amsart} 
\usepackage{amsmath}
\usepackage{amssymb} 
\usepackage{amscd} 
\usepackage{graphicx}  
\usepackage{epsfig,comment}
\usepackage[all]{xy} 
\usepackage{marvosym}
\usepackage{float}
\usepackage{xcolor}
\usepackage{ascmac}
\usepackage[margin=1.2in]{geometry} 
\usepackage{hyperref}

\hypersetup{
	colorlinks=true,
	linktoc=all,
    	linkcolor={red!50!black},
   	citecolor={blue!50!black},
  	urlcolor={blue!80!black}
	}
\usepackage{booktabs} % Enhances quality of tables
\usepackage{tikz-cd,mathtools}

\usetikzlibrary{arrows}
\usetikzlibrary{braids}

\usepackage{wrapfig}
\usepackage{pdfpages}
\usepackage{mwe}

\usepackage{caption}

\DeclareMathOperator{\R}{\mathbb R}

\DeclareMathOperator{\Z}{\mathbb Z}

\DeclareMathOperator{\brl}{brloc}

\DeclareMathOperator{\dc}{dc}

\newcommand{\mc}{\mathcal}

\theoremstyle{plain}
\newtheorem{theorem}{Theorem}
\newtheorem{lemma}[theorem]{Lemma}
\newtheorem{proposition}[theorem]{Proposition}
\newtheorem{corollary}[theorem]{Corollary}
\newtheorem{conjecture}[theorem]{Conjecture}

\newtheorem{question}[theorem]{Question}
\newtheorem{fact}[theorem]{Fact}
\theoremstyle{definition}
\newtheorem{remark}[theorem]{Remark}

\newtheorem{definition}[theorem]{Definition}
\newtheorem{notazione}[theorem]{Notation}
\newtheorem{example}[theorem]{Example}
\newtheorem{exercise}[]{Exercise}
\newtheorem{problem}[theorem]{Problem}

\newtheorem{vuoto}[theorem]{}

\numberwithin{theorem}{section}

\newcommand{\bt}{\begin{theorem}}
\newcommand{\et}{\end{theorem}}

\newcommand{\bv}{\begin{vuoto}}
\newcommand{\ev}{\end{vuoto}}

\newcommand{\bl}{\begin{lemma}}
\newcommand{\el}{\end{lemma}}

\newcommand{\bd}{\begin{definition}}
\newcommand{\ed}{\end{definition}}

\newcommand{\beq}{\begin{equation}}
\newcommand{\eeq}{\end{equation}}

\newcommand{\bexa}{\begin{example}}
\newcommand{\eexa}{\end{example}}

\newcommand{\bexe}{\begin{exercise}}
\newcommand{\eexe}{\end{exercise}}

\newcommand{\bfact}{\begin{fact}}
\newcommand{\efact}{\end{fact}}

\newcommand{\bprop}{\begin{proposition}}
\newcommand{\eprop}{\end{proposition}}

\newcommand{\bp}{\begin{proof}}
\newcommand{\ep}{\end{proof}}

\newcommand{\bc}{\begin{corollary}}
\newcommand{\ec}{\end{corollary}}

\newcommand{\bq}{\begin{question}}
\newcommand{\eq}{\end{question}}

\newcommand{\bcong}{\begin{conjecture}}
\newcommand{\econg}{\end{conjecture}}

\newcommand{\bproblem}{\begin{problem}}
\newcommand{\eproblem}{\end{problem}}

\newcommand{\bs}{\begin{proof}[Proof.]}
\newcommand{\es}{\end{proof}}

\newcommand{\br}{\begin{remark}}
\newcommand{\er}{\end{remark}}

\newcommand{\bn}{\begin{notazione}}
\newcommand{\en}{\end{notazione}}

\begin{document}
\baselineskip 13pt

\title{A formula for the Euler class of foliations}
\author{Alessandro V. Cigna}


\begin{abstract}
    \noindent Given a cooriented branched surface $\mathcal B$ fully carrying a foliation $\mathcal F$, we use the dual graph of $\mathcal B$ to define a simplicial 1-cycle $\Gamma_m(\mathcal B)$ representing the Poincar\'e dual of the Euler class of $\mathcal F$ relative to the boundary. As an example, we complete the classification of which homology classes in the Whitehead link exterior are realisable as relative Euler classes of taut foliations. 
    We also show how our formula generalises previous results of Lackenby and Dunfield. Finally, we observe that cooriented branched surfaces whose complement is a union of balls satisfy a Combinatorial Transverse Surface Theorem, in the sense of Landry--Minsky--Taylor. 
\end{abstract}

\maketitle

\section{Introduction}
In the study of taut foliations on 3-manifolds, the Euler class plays a central role. This paper exhibits a way to compute, understand and manipulate the Euler class of any cooriented foliation, given a branched surface that carries it. 

Given a cooriented foliation $\mc F$ on a compact $3$-manifold $M$, the \emph{Euler class} $e(\mc F)\in H^2(M;\Z)$ is the obstruction to finding a nonvanishing section for the plane bundle $T\mc F\to M$. When $M$ has toroidal boundary and $\mc F$ is transverse to it, we can fix a nonvanishing section $s$ of $TF|_{\partial M}$ and define a relative Euler class $e(\mc F,s)\in H^2(M,\partial M;\Z)$. This time, $e(\mc F,s)$ is the obstruction to extending $s$ to a nonvanishing section over $M$. When the section $s$ is not explicitly defined, we will just write $e(\mc F)=e(\mc F,s)\in H^2(M,\partial M;\Z)$ and we will understand that $s$ is an outward-pointing section.

In \cite{Thurston1986ANF}, Thurston showed that the relative Euler class of a taut foliation is related to the topological complexity of the surfaces sitting inside $M$:

\bt[Thurston's Inequality]\label{thm: thurston's inequality} Let $\mc F$ be a taut foliation of a compact oriented $3$-manifold $M$, with $\mc F$ transverse to $\partial M$. Let $S\subset M$ be a properly embedded oriented surface without sphere or disk components. Then the following inequality holds:
\begin{equation*}
    |\langle e(\mc F), [S]\rangle| \le |\chi(S)|
\end{equation*}
with $e(\mc F)\in H^2(M,\partial M;\Z)$ and $[S]\in H_2(M,\partial M;\Z)$. Moreover, equality holds if $S$ is a union of coherently oriented leaves of $\mc F$.
\et

Another way to think of Thurston's Inequality is through the \emph{Thurston norm}. The Thurston norm of a homology class $\alpha\in H_2(M,\partial M;\Z)$ is the minimal realisable value of $|\chi(S')|$ among all properly embedded orientable surfaces $S$ representing $\alpha$, and $S'$ is $S$ minus its sphere or disk components \cite{Thurston1986ANF}. Thurston's Inequality implies that the relative Euler class of a taut foliation lies in the unit ball of the dual Thurston norm on $H^2(M,\partial M;\Z)$. If the taut foliation has a compact leaf of negative Euler characteristic, then its Euler class actually lies on the boundary of the dual Thurston ball.

The \emph{Euler class-one Conjecture}, posed by Thurston in \cite{Thurston1986ANF}, stated that every class in $H^2(M,\partial M;\Z)$ with unit dual Thurston norm and subject to a necessary parity condition is the relative Euler class of some taut foliation on $M$. The conjecture was proved for vertices of the unit ball of the dual Thurston norm by Gabai \cite{Gabai1997, gabai_fully_2020}. Then the conjecture was disproved by Gabai and Yazdi in \cite{yazdi_thurstons_2020} and \cite{gabai_fully_2020}. By means of the \emph{Fully Marked Surface Theorem}, Yazdi exhibited infinitely many examples of closed hyperbolic $3$-manifolds for which the conjecture does not hold. Lately, Liu showed that counterexamples to the conjecture are not sporadic. Indeed, every oriented closed hyperbolic $3$-manifold $N$ admits some finite cover $M$ for which the conjecture does not hold \cite{liu}. However, we still lack a complete picture of which homology classes are (relative) Euler classes of taut foliations. What other special properties do such classes have?

The Thurston norm is also a valuable tool in its own right. The shape of its unit ball has a close relationship with many other remarkable objects associated with the manifold: fibrations, pseudo-Anosov flows, representations, and Floer homologies are just some of them. Consult \cite{Kitayama2022} for an overview on the Thurston norm and some of its applications. Although this norm is notoriously difficult to compute, Thurston's Inequality forces the unit ball to lie between the two hyperplanes $\{\alpha\in H_2(M,\partial M;\R)\, |\; \langle e(\mc F),\alpha\rangle=\pm 1\}$. 

\bigskip

Understanding the Euler class of a foliation therefore becomes a key problem, both from a conceptual and a computational point of view. The contribution of this paper is a formula for the Euler class of a foliation in terms of a branched surface that carries it.

If $\mc B$ is a cooriented branched surface whose exterior consists of product balls, we endow the dual graph $\Gamma(\mc B)$ of $\mc B$ with a weights structure on its edges and refer to this data as the \emph{maw dual graph} $\Gamma_m(\mc B)$. Theorem \ref{thm: maw dual graph} and Lemma \ref{lemma: graph is cycle} together imply the following result.

\bt\label{thm: maw dual graph intro} Let $\mathcal{B}$ be a cooriented branched surface in $M$ whose exterior consists of product balls. Let $\mathcal{F}$ be a foliation carried by $\mathcal{B}$. The maw dual graph $\Gamma_m(\mc B)$ represents in $H_1(M;\Z)$ the Poincar\'e dual of the relative Euler class $e(\mathcal{F})\in H^2(M,\partial M;\Z)$.
\et

To exemplify a use of this construction, in Section \ref{sec: whitehead link} we define a branched surface $\mc B$ in the exterior of the Whitehead link $L$. We then show that $\mc B$ carries a taut foliation with relative Euler class dual to the oriented meridian of a component of $L$. By putting this together with the works of Gabai \cite{Gabai1997, gabai_fully_2020} and of Fan--Lai--Yu \cite{fan2025noteeulerclass0}, we conclude:

\bigskip

\noindent\textbf{Theorem \ref{thm: whitehead}.} \textit{The relative Euler classes of taut foliations in the Whitehead link exterior are exactly the integral points of unit dual Thurston norm in $H^2(M_L,\partial M_L;\R)$.}

\bigskip

The maw dual graph construction draws attention to cooriented branched surfaces whose exterior is a union of product balls. How much information do such structures capture? By taking inspiration from the work of Landry and Landry--Minsky--Taylor on veering triangulation \cite{LANDRY, LandryMinskyTaylor2024, LandryMinskyTaylor+2026+203+257}, we observe that their Combinatorial Transverse Surface Theorem easily extends to this setting.

\bigskip

\noindent\textbf{Theorem \ref{thm: TST}.}\textit{ Let $M$ be a compact oriented $3$-manifold and $\mc B\subset M$ a cooriented branched surface with generic branched locus and complement a union of balls. An integral homology class $\alpha\in H_2(M,\partial M;\Z)$ pairs nonnegatively with every directed cycle in $\Gamma(\mc B)$ if and only if a surface representative $S$ of $\alpha$ is carried by $\mc B$. Moreover, $S$ has no null-homologous components and, if $\mc B$ carries a taut foliation $\mc F$, $S$  realises the Thurston norm of $\alpha$.}

\bigskip

Sometimes it is useful to understand how the Euler class of a foliation changes when the foliation is carried by a branched surface and the orientation of one of its sectors is reversed to produce a new foliation. A similar situation arises, for instance, from different choices of orientation for a decomposing surface at the end of a taut sutured manifold hierarchy \cite[Lemma 5.6]{cigna2026suturedmanifoldhierarchiesthurston}. The following lemma addresses this situation.

\bigskip

\noindent\textbf{Lemma \ref{lemma: swapping orientation}.} \textit{Let $\mathcal{B'}\subset M$ be a cooriented branched surface with exterior $(N,\gamma)$. Let $(D,\partial D)\subset (N,\partial N)$ be an oriented disk. Consider the cooriented branched surface $\mc B_+\subset M$ obtained by adding $D$ to $\mc B'$ (i.e., by suitably smoothing $D$ on $\mc B'$ about its boundary). Similarly, define the cooriented branched surface $\mc B_-$ by adding $\overline D$ to $\mc B'$. If the complement of $\mc B_
+$ in $M$ is a union of balls, then the symplicial $1$-chain $\Gamma_m(\mc B_+)-\Gamma_m(\mc B_-)$ is homologous to the chain and $(2-|D\cap \gamma|)a_{\mc B_+}(D)$, where $a_{\mc B_+}(D)$ is the oriented arc of $\Gamma(\mc B_+)$ dual to the sector $D$.}

\bigskip

As discussed above, Euler classes of taut foliations are closely related to the Thurston norm. In \cite{cigna2026suturedmanifoldhierarchiesthurston}, we incorporate the maw dual graph construction into an approach for computing the Thurston norm of $3$-manifolds. Forthcoming work of the author produces new counterexamples to the Euler class–one Conjecture for manifolds with toroidal boundary. For example, one can completely classify which classes of dual Thurston norm one arise as relative Euler classes of taut foliations in the exterior of positive $2$-bridge links. In this case as well, the maw dual graph proves a handy tool.

The Euler class is also useful for studying the left-orderability of the fundamental group of a rational homology sphere. Indeed, if such a manifold admits a taut foliation with vanishing Euler class, its fundamental group is left-orderable \cite[Theorem 8.1]{left-order}. Although the converse is not true in general  \cite{boyer20253manifoldsadmittingcoorientabletaut}, it would be of interest to understand when the maw dual graph represents the zero class in homology.

\bq Under which conditions on $\mc B$ and $M$ does the maw dual graph $\Gamma_m(\mc B)$ vanish in homology? 
\eq

The idea of defining a graph dual to the Euler class of a foliation is not new. Lackenby used a graph dual to an ideal triangulation to mimic the role of the Euler class of a taut foliation \cite{Lackenby2000}. In \cite{dunfield_floer_2020}, given a foliar orientation on a one-vertex triangulation of a $3$-manifold, Dunfield represented the Euler class of the corresponding taut foliation by a cocycle supported on a graph. In Section \ref{sec: previous work}, the relation between Lackenby and Dunfield's constructions and the maw dual graph is investigated.

\bigskip

\noindent \textbf{Structure of the paper.} Sections \ref{sec: branched surfaces} and \ref{sec: euler classes} are devoted to recalling the definitions and some of the main properties related to sutured manifolds, branched surfaces and Euler classes. In Section \ref{sec: maw dual graph} we define the maw dual graph and we prove Theorem \ref{thm: maw dual graph} and therefore Theorem \ref{thm: maw dual graph intro}. Then we apply it to establish Lemma \ref{lemma: swapping orientation}. In Section \ref{sec: whitehead link}, we focus on the Whitehead link exterior and prove Theorem \ref{thm: whitehead}. In Section \ref{sec: previous work}, we study the relation between the maw dual graph and previous work linking dual graphs and taut foliations. Finally, in Section \ref{sec: carried vs fully marked}, we inspect which topological information are retained by cooriented branched surfaces with complement a union of balls. There, we prove Theorem \ref{thm: TST}.  

\bigskip

\noindent \textbf{Notation.}  If $P$ is a manifold, $|P|$ is the number of connected components of $P$. If $P$ is oriented and connected, $\overline P$ indicates $P$ with reversed orientation.
The symbol $M$ generically indicates a smooth, compact and oriented $3$-manifold.
If $L$ is a link in $S^3$, the link exterior $S^3-\overset{\circ}{N}(L)$ is referred to as $M_L$.

If $S$ is a properly embedded surface in $M$, we use $N(S)$ to refer to a closed small tubular neighbourhood $S\times[-1,1]$ of $S$ in $M$. The symbol $\partial N(S)$ will indicate $S\times\{\pm 1\}$, whilst $\overset{\circ}{N}(S)$ will indicate $S\times (-1,1)$. The manifold obtained by \emph{cutting} $M$ along $S$ is $M-\overset{\circ}{N}(S)$.

We will understand the integral coefficients for the homology: for example, we will refer to $H_2(M,\partial M;\Z)$ just as $H_2(M,\partial M)$.

\bigskip

\noindent \textbf{Acknowledgements.} I would like to express my gratitude to my advisor, Mehdi Yazdi, for his insightful discussions on the subject of this article and for his constant support. Thanks also to Diego Santoro and Michael Landry for stimulating conversations.

\section{Taut foliations and branched surfaces}\label{sec: branched surfaces}
Taut foliations play a fundamental role in 3–manifold topology, as their existence imposes strong global constraints on the ambient manifold. In particular, a 3–manifold admitting a taut foliation is either $S^2\times S^1$ or is irreducible and has universal cover homeomorphic to $\R^3$. Taut foliations are also closely connected to Thurston’s norm, tight contact structures, and modern gauge-theoretic invariants. For these reasons, understanding when and how taut foliations exist is a central problem.

We now introduce taut foliations and branched surfaces, which will provide the framework for the constructions developed in this paper. 

\bd[Taut foliations] A \emph{foliation} on a compact manifold $M$ is a family $\mathcal F$ of disjoint injectively immersed surfaces - called \emph{leaves} of $\mathcal F$ - that altogether cover $M$ and such that a small neighbourhood of each point is foliated as an open subset of $\R^2\times [0,+\infty)$ subdivided into vertical half-planes $\R\times \{y\}\times[0,+\infty)$.

The foliation $\mathcal F$ is \emph{taut} if it is cooriented and every leaf intersects a closed transversal.
\ed

\begin{figure}
    \centering
    \includegraphics[width=\linewidth]{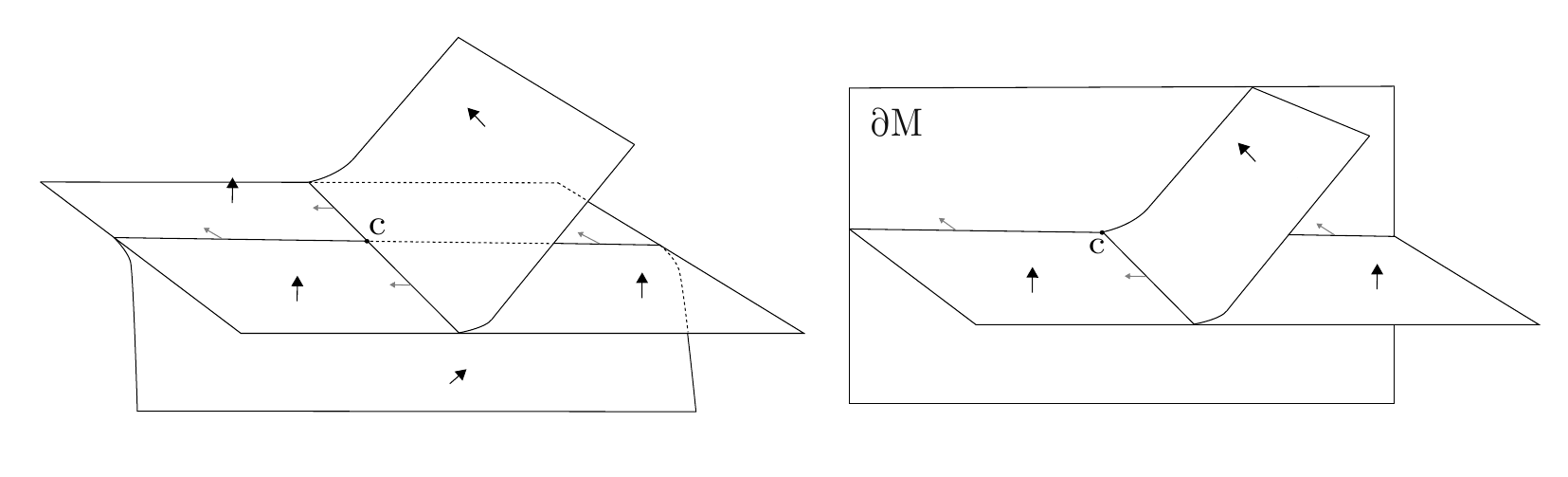}
    \caption{Local views of a branched surface close to a triple point $c$ (left) and at the intersection with $\partial M$ (right). Big arrows indicate a coherent choice of coorientations for the sectors. The small grey arrows indicate the maw vector field. In this case, the normal orientation of $\partial M$ is pointing out of $M$.}
    \label{fig: branching models}
\end{figure}

One classical way to construct foliations on a $3$-manifold $M$ is through branched surfaces. The process goes as follows: define a branched surface $\mathcal B$ on $M$ and show it fully carries a lamination. If the exterior of $\mathcal B$ in $M$ is a product, the lamination can be extended to a foliation. Moreover, by adding hypotheses on $\mathcal B$ (e.g., cooriented, taut, etc.), one can ensure the constructed foliations have good properties. We now present this construction in more detail.

\bd[Branched surface] A \emph{branched surface} in a $3$-manifold $M$ is a $2$-complex $\mathcal B$ that locally looks like an open subset of the models in Figure \ref{fig: branching models}. Notice that a branched surface is naturally endowed with a tangent space at every point. Moreover, we always assume the branched surface to be transverse to $\partial M$.  

The \emph{branching locus} of $\mathcal{B}$ is the set $\brl(\mathcal B)$ of nonmanifold points of $\mathcal B$. It is a properly embedded graph on $\mathcal B$ whose vertices have degree either one or four, and the degree-four vertices are called \emph{triple points}. The branching locus $\brl(\mathcal B)$ is endowed with the \emph{maw vector field}: this is any nonvanishing smooth vector field tangent to $\mathcal B$ and always pointing in the branching direction. See again Figure \ref{fig: branching models}. Once a normal orientation on $\partial M$ is chosen, we extend the maw vector field on $\partial \mc B$ by requiring it always points in the normal direction.

A \emph{sector} of $\mathcal B$ is the metric closure of any connected component of $\mathcal B -\brl(\mathcal B)$.

The branched surface $\mathcal B$ is \emph{cooriented} if there is a choice of normal direction on each sector of $\mathcal B$ that is coherent at the branching locus. Equivalently, $\mathcal B$ is cooriented if its tangent bundle $T\mathcal B$ is.
\ed

\begin{figure}
    \centering
    \includegraphics[width=0.7\linewidth]{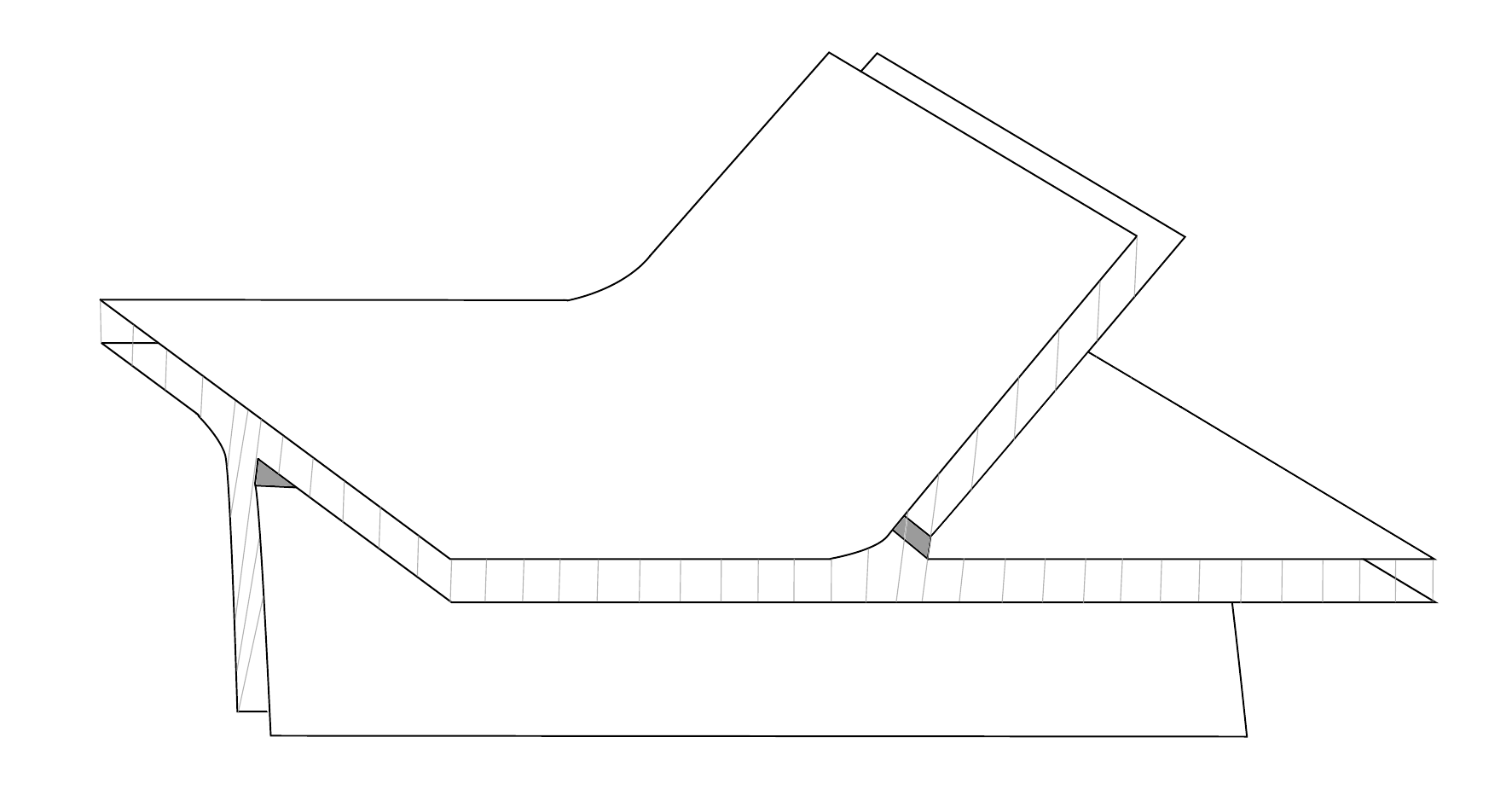}
    \caption{Local view of the fibered neighbourhood $N(\mc B)$ of a branched surface $\mc B$ close to a triple point of the branching locus. The grey intervals are some fibers. The small shaded areas are part of the vertical boundary $\partial_vN(\mc B)$.}
    \label{fig: fibered neighborhood}
\end{figure}

\bd[$I$-fibered neighbourhood and exterior] Let $\mc B$ be a branched surface in a compact $3$-manifold $M$. An \emph{$I$-fibered neighbourhood} of $B$ in $M$ is a regular neighbourhood $N(\mc B)$ foliated by interval fibers that intersect $\mc B$ transversely, see Figure \ref{fig: fibered neighborhood}. We will think of $\mathcal B$ as embedded inside $N(\mathcal B)$, and endowed with a fiber-collapsing map $\pi_{\mc B}:N(\mc B)\to \mc B$. 

The boundary of $N(\mc B)$ can be decomposed into a \emph{horizontal} and a \emph{vertical} part. The \emph{horizontal boundary} $\partial_hN(\mc B)$ is given by the points on $\partial N(\mathcal B)$ that are transverse to the $I$-fibers. If $\mc B$ is cooriented, each component of $\partial_hN(\mc B)$ automatically inherits a normal orientation. The \emph{vertical boundary} $\partial_vN(\mathcal B)$ consists of the portion of $\partial N(\mathcal B)\cap \overset{\circ}M$ that is tangent to the interval fibers.

The \emph{exterior} of $\mathcal B$ in $M$ is the data $(M(\mathcal B),\gamma(\mathcal B))$ where:
\begin{itemize}
    \item $M(\mathcal{B})$ is the closure of $M-N(\mathcal B)$;
    \item $\gamma(\mathcal B)$ is the closure of $\partial M(\mc B)-\partial_hN(\mc B)$. In other words, $\gamma(\mc B)$ is the union of $\partial_vN(\mathcal B)$ and the closure of $\partial M- N(\mathcal B)$. 
\end{itemize}
The connected components of $M(\mathcal B)$ are called \emph{complementary regions} of $\mathcal B$.
We say that a complementary region $N$ of $\mc B$ is a \emph{product} if there exists a compact oriented surface $S$ so that the couple $(N,\gamma(\mc B)\cap N)$ is homeomorphic to $(S\times [0,1],\partial S\times [0,1])$. 
\ed

\bd[Laminations] A \emph{lamination} in a compact $3$-manifold $M$ is a family of disjoint injectively immersed surfaces which altogether foliate a closed subset of $M$.

Given a branched surface $\mathcal B\subset M$, we say that a lamination $\mathcal{L}$ is \emph{carried} by $\mathcal B$ if it can be isotoped to lie inside an $I$-fibered neighbourhood $N(\mathcal B)$ in a way that each leaf of $\mathcal L$ is transverse to the $I$-fibers. We say that $\mathcal L$ is \emph{fully carried} by $\mathcal B$ if, furthermore, each interval fiber intersects at least one leaf of $\mathcal L$.
\ed

As already anticipated, branched surfaces can be used to construct foliations. Indeed, suppose that $\mathcal B\subset M$ is a branched surface that fully carries a lamination $\mc L$. Then, after possibly doubling some leaves, we can isotope $\mc L$ so that the horizontal boundary $\partial _hN(\mc B)$ is contained in a union of leaves of $\mc L$. Since $\mc L$ is fully carried by $\mc B$, when we cut $N(\mc B)$ along $\mc L$, we get an interval-bundle over a surface. In particular, $\mc L$ can be extended to a foliation $\mc F_0$ of $N(\mc B)$ transverse to the interval fibers, and such that $\mc F_0$ restricts to a product foliation on the vertical boundary $\partial_vN(\mc B)$. If the exterior of $\mc B$ is a union of product components, we can further extend $\mathcal F_0$ to a foliation $\mc F$ on $M$ through the product foliations of the complementary regions. We will briefly say that $\mc B$ fully carries (or just carries) the foliation $\mc F$.

If $\mc B$ is cooriented, then $\mathcal F$ is cooriented too. If every sector of $\mc B$ intersects a closed positive transversal, then the foliation $\mc F$ is taut.

\bigskip

The construction above motivates the importance of determining whether a branched surface fully carries a lamination. A classical example is when the branching locus of $\mc B$ has no triple points. In this case, a fully carried lamination can be constructed in the following way. For every sector $s$ of $\mc B$, embed a Cantor set $C_s$ of parallel copies of $s$ inside $N(\mc B)$. More precisely, $C_s$ is homeomorphic to $s\times C$, where $C$ is the Cantor set, and $C_s$ is embedded in $\pi_{\mc B}^{-1}(s)\subset N(\mc B)$. By the self-similarity of $C$, the surfaces $\bigcup_{\text{sectors } s}C_s$ can be glued together along their boundaries to give rise to a lamination fully carried by $\mc B$. 

When the branching locus of $\mc B$ has triple points, it can happen that $\mc B$ fully carries no lamination. We refer the interested reader to \cite{li_laminar_2002} for a very general sufficient criterion to verify whether a generic branched surface fully carries a lamination. 

\br Given a coorientable foliation $\mc F$, one can always construct a coorientable branched surface that carries $\mc F$. This is the content of Christy and Goodman's construction as explained in \cite[Section I]{pjm/1102364191}. The idea is to fix a transverse flow $\phi$ to $\mc F$ and a collection of compact planar surfaces $\Delta=\{D_i\}_i$ contained in the leaves of $\mc F$. If the collection $\Delta$ is chosen suitably (e.g., every orbit of $\phi$ needs to intersect the interior of some surface in $\Delta$), a branched surface $\mc B$ is obtained by collapsing the restricted orbits of $\phi$ in $M-\overset{\circ}{\Delta}$. Moreover, the collection $\Delta$ can be chosen to consist of disks, in which case the complement of $\mc B$ is a union of balls.
\er

\section{Euler classes of plane bundles}\label{sec: euler classes}

Given an oriented circle bundle $\pi:E\to X$ on a cellular complex, the Euler class $e(\pi)$ is an element of $H^2(X)$ measuring the obstruction to finding a section of $\pi$, hence how far $\pi$ is from being the trivial circle bundle. We will now define the Euler class and state the main properties that we will need throughout the paper. A good resource to learn the details of these results is \cite[Chapter 4]{candel_foliations_2003}, which the reader is suggested to consult in case of further interest.

\bigskip

In what follows, it is useful to bear in mind the following principle. If $S$ is a circle or a compact oriented surface with non-empty boundary, then its unit tangent bundle $T_1S$ is trivial:  there exists a section $S\to T_1S$ or, equivalently, there exists an isomorphism of circle-bundles $\psi:T_1S\to S\times S^1$. If two such trivialisations are given, say $\psi_1$ and $\psi_2$, then we can compare them via a map $\phi:S\to S^1$ defined by the following composition: $$S\xhookrightarrow{i_S}S\times S^1\overset{\psi_2^{-1}}\longrightarrow T_1S\overset{\psi_1}\longrightarrow S\times S^1\overset{\pi_{S^1}}\longrightarrow S^1.$$ We will refer to $\phi$ as the map $S\to S^1$ \emph{induced from $\psi_2$ with respect to $\psi_1$}.

\bigskip

Let $\pi:E\to X$ be an oriented circle bundle over a finite $n$-dimensional cellular complex $X$. Let $X^i$ be the $i$-th skeleton of $X$. We attempt to define a section $s$ to $\pi$ by proceeding inductively on the skeletons. Choose arbitrarily a section $s_0$ on the vertices $X^0$, then extend it continuously to a section $s_1$ on $X^1$. Given a $2$-cell $\Delta\subset X^2$, we would like to extend the section $s_1$ over $\Delta$. Since $\Delta$ is contractible, the restriction of $\pi$ to it is trivial, and therefore isomorphic to $\Delta\times S^1$. Through this identification, the section $s_1$ on the boundary induces a map $\partial \Delta\to S^1$ with degree $c_{s_1}(\Delta)$. Once $c_{s_1}(\Delta)$ is defined for every $2$-cell $\Delta$, we have a cellular $2$-cochain $c_{s_1}$ in $C^2_{CW}(X)$.

The following properties are classical. For a proof, consult for example \cite[Section 4.3]{candel_foliations_2003}.

\bt\label{thm: obstruction theory} The cochain $c_{s_1}$ satisfies: \begin{enumerate}
    \item $c_{s_1}$ does not depend on the choice of the trivializations $E_{|\Delta}\cong \Delta\times S^1$.
    \item $c_{s_1}$ is in fact a $2$-cocycle, hence representing a class $[c_{s_1}]\in H^2(X)$.
    \item If $t_1$ is another section of $\pi:E_{|X^1}\to X^1$, then $c_{s_1}$ and $c_{t_1}$ are cohomologous.
    \item The section $s_1$ can be extended to a section of $E_{|X^2}$ if and only if $[c_{s_1}]=0\in H^2(X)$, if and only if it can be extended to a section of $\pi:E\to X$.
\end{enumerate} 
\et

Suppose now that we know a section $\sigma:Y\to E_{|Y}$ defined on a subcomplex $Y$ of $X$. We might wonder whether $\sigma$ can be extended to a section over the whole $X$ and proceed in a similar fashion as above. Namely, pick an arbitrary section $s_1$ on $X^1\cup Y$ that coincides with $\sigma$ on $Y$ and define for every $2$-cell $\Delta$ in $X-Y$ the number $c_{s_1}(\Delta)$. As we are fixing the section on $Y$ and we just care about cells not in $Y$, the integer-valued function $c_{s_1}$ gives this time a cellular relative cochain in $C^2_{CW}(X,Y)$. The analogue of Theorem \ref{thm: obstruction theory} holds also in this more general context, guaranteeing that $c_{s_1}$ induces a relative cocycle in $H^2(X,Y)$ that vanishes if and only if $\sigma$ can be extended to a section over $X$.

\bigskip

\subsection{The Poincar\'e--Hopf Theorem}

Let $S$ be a compact oriented surface, possibly with boundary, and endowed with a Riemannian metric. Consider the unit tangent bundle $\pi:T_1S\to S$ and a section $m$ defined on $\partial S$. The Poincar\'e--Hopf Theorem establishes that one can compute the Euler class of $\pi$ by drawing a suitable vector field on $S$. 

\bd[Index] Let $X$ be a vector field on $S$ possibly with isolated zeroes in the interior of $S$. If $z\in S$ is a zero of $X$, then pick a small disk $D_z\subset S$ having only the zero $z$ in it. Fix a trivialisation of $\pi$ over $D_z$ and consider the map $X_{|\partial D_z}:\partial D_z\to S^1$ induced by the section $X|_{\partial D_z}$. The \emph{index} of $X$ at $z$ is the degree of $X|_{\partial D_z}$. Call \emph{index sum} of $X$ the number $\iota_X$ given by the sum of all the indices of $X$ at its zeroes. 
\ed

\bl\label{lemma: degree bordism} Let $S$ be a compact oriented surface with non-empty boundary and $m$ a nonvanishing vector field on $\partial s$. Fix a trivialisation for $\pi:T_1S\to S$ and call $c_m(S)$ the degree of the map $m_{|\partial S}:\partial S\to S^1$ induced by the section $m$. If $X$ is a smooth vector field on $S$ with isolated zeroes in $S$ and coinciding with $m$ at the boundary, then $$c_m(S)=\iota_X.$$
\el
\bp As in the previous definition, pick a small disk $D_z$ for every singularity $z$ of $X$. After possibly shrinking, assume that the disks $D_z$ are pairwise disjoint. The restriction of $X$ to $B:=S-\cup_z\overset{\circ}{D_z}$ is nonvanishing, hence induces another trivialisation of the unit tangent bundle of $B$. In particular, we obtain a map $\phi: B\to S^1$ and a homomorphism $\phi_*: H_1(B)\to H_1(S^1)$. Since $\overline {\partial S}$ and $\cup_z \partial D_z$ are homologous in $B$, their images through $\phi_*$ coincide. But the integer $\phi_*([\partial S])\in H_1(S^1)\cong \Z$ equals the degree $c_m(S)$ and similarly $\phi_*[\cup_z \partial D_z]$ equals $-\iota_X$. The lemma follows.
\ep

The following is a very well-known result in differential topology and several proofs are available in the literature. One such source is \cite[Theorem 14.4]{benedetti_lectures_2021}.

\bt[Poincar\'e--Hopf Theorem] Let $S$ be a compact oriented surface and let $X$ be a smooth vector field on $S$ transverse to the boundary. If $X$ has isolated zeroes, then $$\langle e(\pi),[S]\rangle=\iota_X=\chi(S),$$ where $e(\pi)\in H^2(S,\partial S)$ is the Euler class of $\pi: T_1(S)\to S$ relative to any boundary section transverse to the boundary, and $[S]\in H_2(S,\partial S)$ is the fundamental class.
\et

\bigskip

\subsection{Euler classes of foliations} 

\bd[Euler class of a foliation] Let $M$ be a compact oriented manifold with toroidal boundary, and let $\mathcal F$ be a foliation on $M$ transverse to $\partial M$. Fix a cellular structure on $M$ and consider the tangent bundle $T\mathcal F$ on $M$. By means of a Riemannian metric for $T\mathcal F\to M$, we can study the unit circle bundle $\pi:T_1\mc F\to M$. Since $\mathcal F$ is transverse to $\partial M$, $\pi|_{\partial M}$ admits an outward pointing section $\sigma$. The \emph{Euler class} of $\mc F$ is the relative $2$-cocycle $$e(\mc F)=[c_{s_1}]\in H^2(M,\partial M)$$ where $s_1$ is any section of $\pi$ on $M^1\cup \partial M$ coinciding with $\sigma$ on $\partial M$. In other words, $e(\mc F)$ is the Euler class of $\pi$ relative to an outward pointing section on $\partial M$.
\ed

We might as well define the Euler class $e(\mc F,-\sigma)$ relative to the inward pointing section $-\sigma$ on $\partial M$. The vector fields $\sigma$ and $-\sigma$ are homotopic through sections of $\pi_{|\partial M}$. Indeed, they differ by a $180^o$ clockwise rotation on each circle of the oriented bundle $T_1\mc F|_{\partial M}$. Therefore, $$e(\mc F,\sigma)=e(\mc F,-\sigma).$$ 

\section{The maw dual graph}\label{sec: maw dual graph}

Let $M$ be a compact oriented $3$-manifold whose boundary is a (possibly empty) union of tori. Fix a transverse orientation for $\partial M$. Let $\mathcal B$ be a cooriented branched surface in $M$ whose exterior is a union of product balls.  In this section, we construct a $1$-cycle $\Gamma_m(\mathcal{B})\subset M$ transverse to $\mathcal B$. When $\mathcal B$ fully carries a lamination, then such a lamination can be extended to a cooriented foliation $\mathcal{F}$ of $M$, and the $1$-cycle $\Gamma_m(\mathcal{B})\in H_1(M)$ is dual to the relative Euler class $e(\mathcal F)\in H^2(M,\partial M)$. For example, this is the case when $\mathcal{B}$ comes from a groomed manifold hierarchy, or when $\mathcal{B}$ is laminar.

\bigskip

Here is the idea of the construction: since the exterior of $\mathcal B$ is a union of balls, we can define a cellular structure on $M$ having $\mathcal B$ as the $2$-dimensional skeleton. The sectors of $\mathcal{B}$ are a union of $2$-cells, the branching locus and $\partial B$ are part of the $1$-skeleton. If $\mathcal{B}$ carries a foliation $\mathcal{F}$, the tangent plane field $T\mc F|_\mc B$ is homotopic to $T\mc B$. Hence, the maw vector field gives a nonvanishing section of the tangent plane field for most of the $1$-skeleton. 
The cycle $\Gamma_m(\mathcal B)$ is an embedded graph such that, for every sector $s$ of $\mathcal B$, the algebraic intersection number between $\Gamma_m(\mathcal B)$ and $s$ equals the obstruction to extending the maw vector field from $\partial s$ to the interior of $s$. Such obstruction is what we will call the \emph{maw Euler characteristic} of $s$.

\begin{figure}
    \centering
    \includegraphics[width=0.8\linewidth]{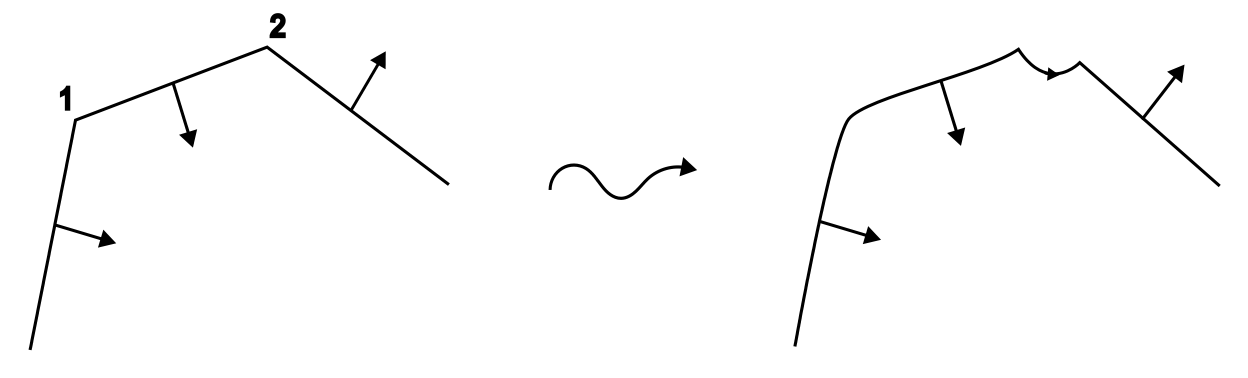}
    \caption{Left: a smoothable vertex, labelled (1), and a double corner, labelled (2). Right: the maw vector field induces a corner structure on $\partial s$: every smoothable vertex can indeed be smoothed (therefore inducing no corners), and each double vertex gives rise to a couple of right angles.}
    \label{fig: corners}
\end{figure}

\bd[Maw Euler characteristic] Let $s$ be a sector of a branched surface $\mathcal B$ and let $v$ be a vertex of the branching locus $\brl (\mc B)$ on $\partial s$, i.e., $v$ is either a triple point on $\partial s$ or an extreme of $s\cap \partial M$. We say that $v$ is a \emph{double corner} (or double vertex) for $s$ if the maw vector field points inward on one edge at $v$ and outward on the other. We say that $v$ is a \emph{smoothable corner} (or smoothable vertex) otherwise. See Figure~\ref{fig: corners} left. Indicate with $\dc(s)$ the number of double corners on $\partial s$. The \emph{maw Euler characteristic} of $s$ is $$\chi_m(s):=\chi(s)-\frac 12 \dc(s).$$
\ed

\begin{figure}
    \centering
    \includegraphics[width=0.8\linewidth]{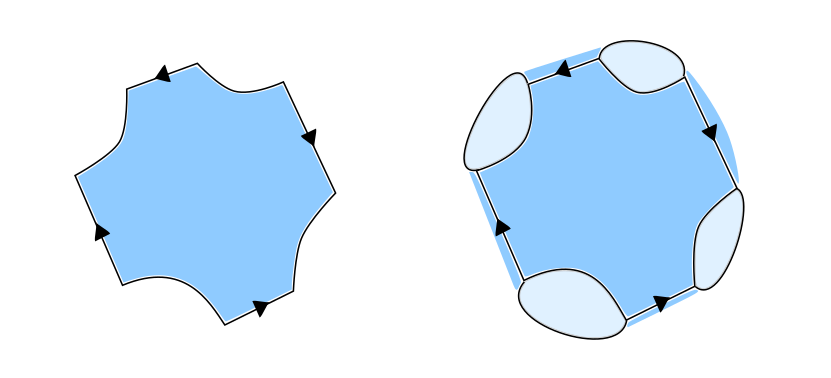}
    \caption{The sector $s$ on the left, and its double $ds$ on the right. The surface $ds$ is obtained by attaching two copies of $s$ along the segments of $\partial s$ where the maw vector field is tangential.}
    \label{fig: doubling}
\end{figure}

\bl\label{lemma: euler class equals maw chi} Let $s$ be a sector of a branched surface $\mathcal B$. Let $m$ be the maw vector field on the boundary of $s$. Suppose $X$ is a smooth vector field with isolated zeroes in $s$ that coincides with $m$ on $\partial s$. Then \begin{equation*}
    \iota_X=\chi_m(s).
\end{equation*}  
\el
\bp Consider $s$ as a surface with corners by introducing two right angles on $\partial s$ for every double vertex as in Figure~\ref{fig: corners} right, where we have slightly isotoped $m$ to be tangent to the new edges of $\partial s$. Apply the same modification to $X$ so that it coincides with the modified $m$ in a neighbourhood of $\partial s$. Call $\partial _\tau s$ the union of the new edges of $\partial s$ which are tangent to $m$. Consider the smooth surface $ds=s\cup\bar s$ obtained by doubling $s$ along $\partial _\tau s$, as in Figure \ref{fig: doubling}. Let $\overline X$ be the vector field on $\bar s$ equal to $-X$ on each point of $\bar s$. The vector fields $X$ and $\overline X$ patch together to a well-defined vector field $dX$ on $ds$ transverse to the boundary. By construction $\iota_X=\iota_{\overline X}$. The Poincar\'e--Hopf Theorem implies $$\chi(ds)=\iota_X+\iota_{\overline X}=2\iota_X.$$
The proof follows by the fact that $\chi(ds)=2\chi(s)-|\partial_\tau s|=2\chi(s)-dc(s).$  
\ep

We can now define the maw dual graph.

\bd[Directed and maw dual graph]\label{def: maw dual graph} Let $\mathcal B$ be a cooriented branched surface in an ambient $3$-manifold $M$ whose complement consists of $3$-balls. Fix a base point on each sector and inside each complementary region of $\mathcal B$. The \emph{directed dual graph} of $\mathcal B$ is the oriented graph $\Gamma(\mc B)$ in $M$ constructed as follows:
\begin{itemize}
    \item Vertices: the vertices of $\Gamma(\mc B)$ are the base points in the complementary regions of $\mathcal B$;
    \item Edges: for every sector $s$ of $\mathcal{B}$, there is an arc $a(s)$ joining the vertices corresponding to the complementary regions of $\mathcal{B}$ adjacent at $s$. The arc $a(s)$ intersects $\mc B$ transversely and only at the base point of $s$. Moreover, $a(s)$ is oriented as the coorientation of $s$.
\end{itemize}

The \emph{maw dual graph} of $\mc B$ is the graph $\Gamma_m(\mc B)$ obtained by adding a system of weights to the edges of $\Gamma(\mc B)$. Namely, for each sector $s$ of $\mc B$, the dual arc $a(s)$ carries a weight equal to $\chi_m(s)$
\ed

\bl\label{lemma: graph is cycle} Let $\mathcal{B}$ be a cooriented branched surface in $M$ whose exterior is a union of product balls. The maw dual graph $\Gamma_m(\mathcal B)$ is a simplicial $1$-cycle. 
\el
\bp Let $R_+$ (resp. $R_-$) be the portion of $\partial_h N(\mc B)$ whose coorientation points inside (resp. outside) $N(\mc B)$. We need to show that at each vertex $v$ of $\Gamma_m(\mathcal B)$, the total weight of the edges pointing into $v$ is the same as the total weight of the edges pointing out of it. The vertex $v$ is the base point of a product ball $B$ in the exterior of $\mathcal B$. The edges pointing out of $v$ correspond altogether with the sectors of $\mathcal B$ that carry $R_+\cap B$, which is a disk. As vector fields on the sectors of $\mathcal{B}$ restricting to the maw vector field at the boundary patch together to a vector field on their union and the index formula is additive, the total weight pointing out of $v$ is $\chi(R_+\cap B)$ (in this case equal to $1$). Similarly, the total weight pointing into $v$ is $\chi(R_-\cap B)$. As $\chi(R_+\cap B)=\chi(R_-\cap B)$ the claim follows.
\ep

\bt\label{thm: maw dual graph} Let $\mathcal{B}\subset M$ be a cooriented branched surface whose exterior is a union of product balls. Let $\mathcal{F}$ be a foliation fully carried by $\mathcal{B}$ and let $\Gamma$ be an oriented weighted graph in $M$. Suppose that $\Gamma$ is a $1$-cycle and intersects each sector $s$ exactly $\chi_m(s)$ times according to orientation (e.g., if $\Gamma=\Gamma_m(\mathcal{B})$). Then $\Gamma$ represents in $H_1(M)$ the Poincar\'e dual of the relative Euler class $e(\mathcal{F})\in H^2(M,\partial M)$.
\et
\bp As the exterior of $\mathcal B$ is a union of product balls, each sector of $\mathcal B$ embeds in the horizontal boundary of some such ball. Therefore, each sector of $\mathcal B$ is a planar surface with boundary. On every sector, choose a finite family of disjoint properly embedded arcs whose complement is a $2$-cell. By adding these arcs to the branching locus $\brl(\mathcal B)$ and to $\partial\mc B$, we obtain the $1$-skeleton of a cellular structure on $M$ having one $2$-cell $s_0$ for every sector $s$ of $\mathcal B$, and one $3$-cell for every complementary ball. Extend the maw vector field to the $1$-skeleton of this cellulation by choosing arbitrarily a vector field transverse to the new arcs and tangent to $\mathcal B$. In order to carry out this extension, we first might need to slightly isotope the maw vector field close to the endpoints of the newly added arcs. Call the extended vector field $m$ again. By definition of the relative Euler class, there is a relative $2$-cochain $c_m$ representing $e(\mc F)$, and such that $c_m(s_0)$ is the obstruction to finding a nonvanishing section of $T\mathcal F$ on $s_0$ that coincides with $m$ on $\partial s_0$. Since $s_0$ is contained in $\mathcal B$, the plane bundle $T\mathcal F|_{s_0}$ is homotopic to the tangent bundle of $s_0$ and so $$c_m(s_0)=\iota_X=\chi_m(s),$$ where $X$ is a smooth vector field with isolated zeroes in $s$ that coincides with $m$ on $\partial s$. The first equality comes from Lemma \ref{lemma: degree bordism}, the second equality is the content of Lemma~\ref{lemma: euler class equals maw chi}. Since $c_m$ and the algebraic intersection cochain $\langle \cdot,\Gamma\rangle$ coincide on the free basis of the cellular $2$-chains of $M$, we conclude that $e(\mc F)=\langle \cdot,[\Gamma]\rangle$.
\ep

\br \begin{itemize}
    \item The proof of Theorem \ref{thm: maw dual graph} extends to more general (non-generic) branched surfaces, provided a maw vector field is well defined. For instance, this is used in the proof of Lemma \ref{lemma: lackenby's graph}, where the branched surface locally looks like Figure \ref{fig: flattening} right.

    \item Theorem \ref{thm: maw dual graph} also extends to the context of sutured manifolds $(M,\gamma)$. In that case, the tangential boundary $R(\gamma)$ is part of the branched surface $\mc B$, whilst $\gamma$ is transverse to it. The maw dual graph $(\Gamma_m(\mc B),\partial \Gamma_m(\mc B))\subset (M,R(\gamma))$ represents the dual in $H_1(M,R(\gamma))$ of the relative Euler class $e(\mc F)\in H^2(M,\gamma)$, as defined in \cite{yazdi_thurstons_2020}. 

    \item The fact that the exterior of $\mc B$ is a union of balls is not enough to guarantee that each sector of $\mc B$ is topologically a disk. A trivial counterexample is given by nesting a trivial bubble on a sector of any branched surface satisfying the assumptions of Theorem \ref{thm: maw dual graph}. A less trivial (although non-taut!) counterexample can be constructed on $S^2\times S^1$. Let $T$ be an oriented torus dividing $S^2\times S^1$ into two solid tori $V_1$ and $V_2$ glued along their meridional curves. Call $D_1$ and $D_2$ two oriented meridional disks in $V_1$ and $V_2$ respectively, so that $\partial D_1\cap \partial D_2= \emptyset$. Let $\mc B$ be the branched surface obtained by smoothing $T\cup D_1\cup D_2$ according to coorientation. The branched surface $\mc B$ fully carries a Reeb foliation of $S^2\times S^1$, the exterior of $\mc B$ is a union of two product balls, yet $\mc B$ has two annular sectors on $T$.    
\end{itemize} 
\er

\subsection{Producing nonisotopic foliations}
Sometimes, we might want to understand how the Euler class of a foliation changes when it is carried by a branched surface, and we reverse the orientation of one of its sectors to construct a different foliation. This operation surely yields a new cooriented branched surface if the given sector is a \emph{source sector}, i.e., the maw vector field points always out of it. 

\bl\label{lemma: swapping orientation} Let $\mathcal{B'}\subset M$ be a cooriented branched surface with exterior $(N,\gamma)$. Let $(D,\partial D)\subset (N,\partial N)$ be an oriented disk. Consider the cooriented branched surface $\mc B_+\subset M$ obtained by adding $D$ to $\mc B'$ (i.e., by suitably smoothing $D$ on $\mc B'$ about its boundary). Similarly, define the cooriented branched surface $\mc B_-$ by adding $\overline D$ to $\mc B'$. If the complement of $\mc B_
+$ in $M$ is a union of balls, then the symplicial $1$-chain $\Gamma_m(\mc B_+)-\Gamma_m(\mc B_-)$ is homologous to the chain and $(2-|D\cap \gamma|)a_{\mc B_+}(D)$, where $a_{\mc B_+}(D)$ is the oriented arc of $\Gamma(\mc B_+)$ dual to the sector $D$.
\el

\begin{figure}
    \centering
    \includegraphics[width=0.7\linewidth]{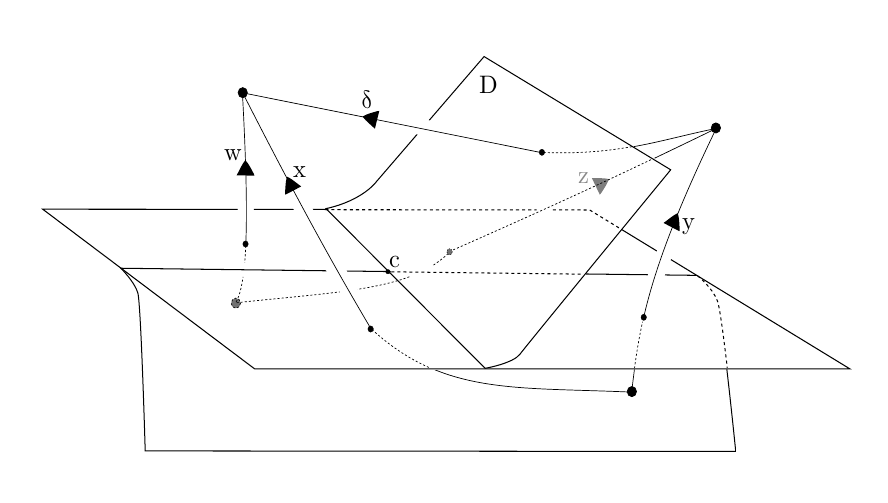}
    \caption{The local view of $\mathcal B_+$ around a triple point $c$ on $\partial D$ but not on $\gamma$. Here, $a_+(c)=x+z$ and $a_-(c)=y+w$. Since $x-y$ and $w-z$ are homologous to $\delta$, the contribution $a_+(c)-a_-(c)$ is zero in homology.}
    \label{fig: four sectors 1}
\end{figure}

\begin{figure}
    \centering
    \includegraphics[width=0.495\linewidth]{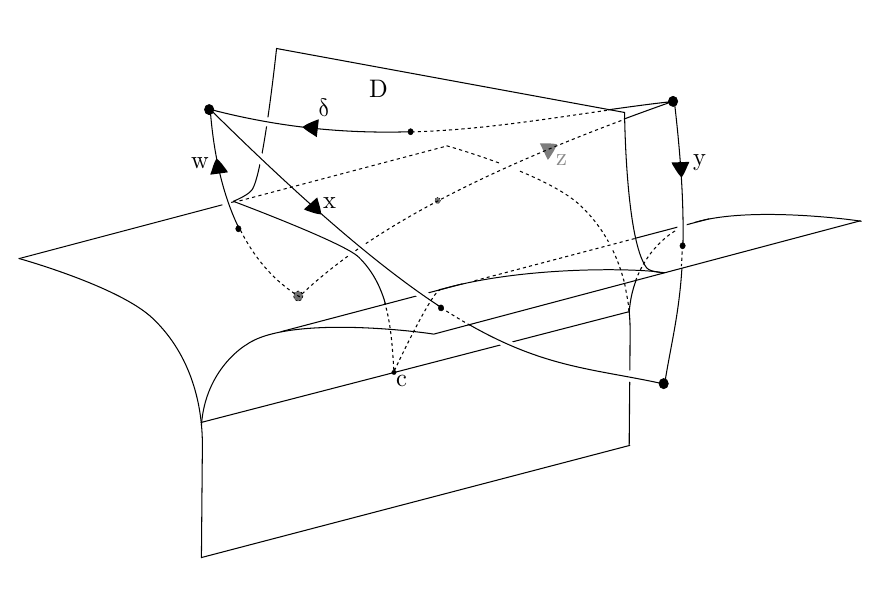}
    \includegraphics[width=0.495\linewidth]{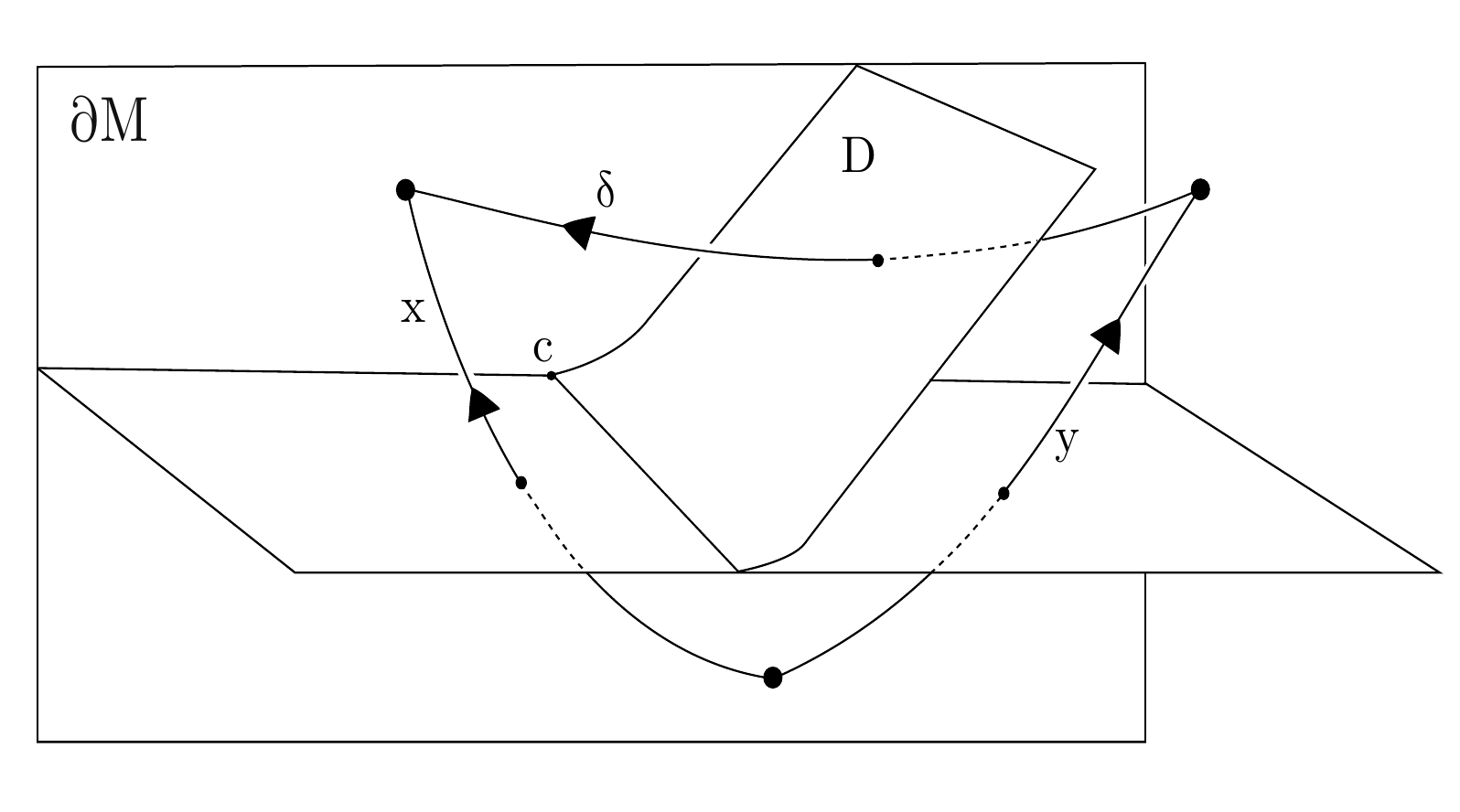}
    \caption{Left: the local view of $\mathcal B_+$ around a triple point $c$ on $\partial D\cap \gamma$. Here, $a_+(c)= y+w$ and $a_-(c)=x+z$. As $y-x$ and $w-z$ are homologous to $\delta$, $a_+(c)-a_-(c)$ is homologous to $2\delta$. Right: the local view of $\mathcal B_+$ around a vertex $c\in D\cap \partial B_+$. Here, $a_+(c)= x$ and $a_-(c)=y$. As $x-y$ is homologous to $\delta$, $a_+(c)-a_-(c)$ is homologous to $\delta$. }
    \label{fig: three sectors}
\end{figure}

\bp Without loss of generality, assume that each component of $\partial M$ has normal orientation pointing out of $M$. As the sectors and complementary components of the two branched surfaces are in one-to-one correspondence, we will refer to them as pairwise identified. Let $\Gamma_{\pm}:=\Gamma_m(\mathcal B_\pm)$ be the maw dual graph constructed as in Definition \ref{def: maw dual graph} after choosing the same family of base points for $\mathcal B_+$ and $\mathcal B_-$. We compute the difference $[\Gamma_m(\mc B_+)-\Gamma_m(\mc B_-)]$ by a closer look at the variation of $\chi_m$ among the sectors. Remember that for every sector $s$ on $\mathcal B_+$ there is a dual arc $a_+(s):=a_{\mc B_+}(s)\subset \Gamma_+$. Call $\delta:=a_+(D)$. In this way, $$[\Gamma_m(\mc B_+)]=\left[ \sum_{s \text{ sector of } \mathcal B_+} \chi_m^{\mathcal B_+}(s) a_+(s)\right]=\left[ \sum_{s \text{ sector of } \mathcal B_+} \chi(s) a_+(s) -\frac 12 \sum_{c \text{ vertex of } \brl(\mathcal B_+)} a_+(c)\right]$$ where, at every vertex $c$ of the branching locus of $\mathcal B_+$, $a_+(c)$ is the sum of the arcs dual to the sectors of $\mathcal B_+$ forming a double corner at $c$. There are exactly two such sectors for every vertex $c$ that is a triple point of $\mathcal B_+$, and one for every vertex $c$ on $\partial \mathcal B_+$. Since an analogous formula holds for $[\Gamma_m(\mc B_-)]$, where also $a_-(s)=a_+(s)$ for every $s\ne D$, $a_-(D)=-a_+(D)=-\delta$, and $a_-(c)=a_+(c)$ for $c\notin \partial D$, we get \begin{equation}\label{eq: difference}[\Gamma_m(\mc B_+)-\Gamma_m(\mc B_-)]=\left[2\delta-\frac 12\sum_{c \text{ vertex on } \partial D} a_+(c)-a_-(c)\right].\end{equation} 
If $c$ is a vertex of $\brl(\mc B_+)$ on $D$, there are three possibilities: 
\begin{itemize}
    \item[$(i)$] $c$ is a triple point on $\partial D$ but not on $\partial D\cap \gamma$ (Figure \ref{fig: four sectors 1}),
    \item[$(ii)$] $c\in\partial D\cap \gamma$ (Figure \ref{fig: three sectors} left), or
    \item[$(iii)$] $c\in D\cap \partial B_+$ (Figure \ref{fig: three sectors} right).
\end{itemize}

     In case $(i)$, the quantity $a_+(c)-a_-(c)$ is an alternate sum of the arcs dual to the four sectors having a vertex at $c$. This quantity is in fact zero in homology, see Figure \ref{fig: four sectors 1}. 

In case $(ii)$, then $a_+(c)-a_-(c)$ is $2\delta$ in homology as in Figure \ref{fig: three sectors} left. 

In case $(iii)$, $a_+(c)-a_-(c)$ is instead homologous to $\delta$, as in Figure \ref{fig: three sectors} right. Since there are two such vertices for every component of $D\cap \partial M$, the lemma follows by Identity (\ref{eq: difference}). 
\ep

\section{Taut foliations in the Whitehead link exterior}\label{sec: whitehead link}

Let $M$ be a compact oriented irreducible $3$-manifold $M$ with toroidal boundary. Which classes in $H^2(M,\partial M;\R)$ are relative Euler classes of taut foliations of $M$? Two important necessary conditions are due to Thurston. 

\vspace{0.2cm}

\noindent \textbf{Dual Thurston norm condition.} \cite{Thurston1986ANF} \textit{The relative Euler class of a taut foliation has dual Thurston norm at most one. This means that for every homology class $\alpha\in H^2(M,\partial M;\R)$, the following inequality holds $$|\langle e(\mc F),\alpha\rangle|\le x(\alpha),$$ where $x$ is the Thurston norm of $\alpha$.}

\bigskip

\begin{figure}
    \centering
    \includegraphics[width=0.3\linewidth]{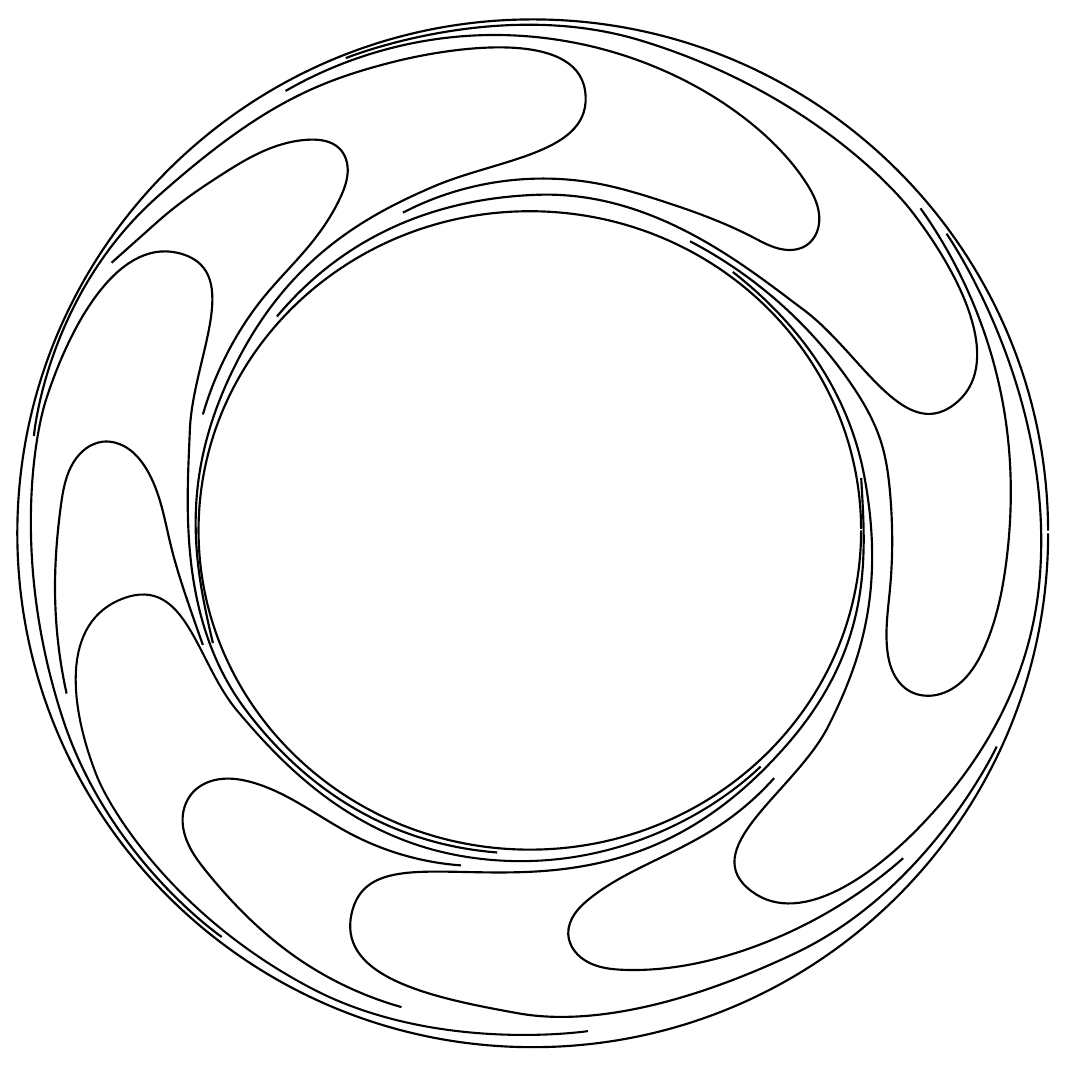}
    \caption{A $2$-dimensional Reeb component is a subset of leaves of a foliation of a surface covering an annulus. The boundary of the annulus is composed of leaves with the same coorientation, i.e., they both point inside or both outside the annulus. The interior leaves are lines.}
    \label{fig: 2-dim Reeb component}
\end{figure}

In general, this is the only necessary condition we know about relative Euler classes of taut foliations. However, if we require the foliation induced on the boundary to have no Reeb annuli (Figure \ref{fig: 2-dim Reeb component}), then the relative Euler class satisfies a necessary \emph{parity condition} as well.

\bigskip

\noindent \textbf{Parity condition.} \cite{Thurston1986ANF} \textit{Let $\mc F$ be a taut foliation on $M$ such that $\mc F\cap \partial M$ has no Reeb components. If $S\subset M$ is a properly embedded oriented surface such that no component of $\partial S$ bounds a disk in $\partial M$, then $$\langle e(\mc F),[S]\rangle \equiv \chi(S) \mod2.$$} 

\bigskip

The Euler class-one Conjecture affirmed that every integral class in $H^2(M,\partial M;\R)$ satisfying the parity condition and with unit dual Thurston norm is the relative Euler class of a taut foliation on $M$ \cite[Conjecture 3]{Thurston1986ANF}. Gabai gave the first partial positive answer to this conjecture:

\bt\cite{Gabai1997, gabai_fully_2020}\label{thm: gabai vertex} Let $M$ be a compact oriented irreducible $3$-manifold, possibly with toroidal boundary, and let $\beta\in H^2(M,\partial M)$ be a vertex of the dual unit ball. Then there is a taut foliation on $M$, without Reeb components at the boundary, whose relative Euler class is $\beta$. 
\et

The conjecture was then definitely proved false by Yazdi and Gabai--Yazdi in \cite{yazdi_thurstons_2020} and \cite{gabai_fully_2020}. In particular, the original question still remains unanswered in its generality: which homology classes are relative Euler classes of taut foliations? The answer to this question is not known even in simple cases. 

\bigskip

\begin{figure}
    \centering
    \includegraphics[width=0.4\linewidth]{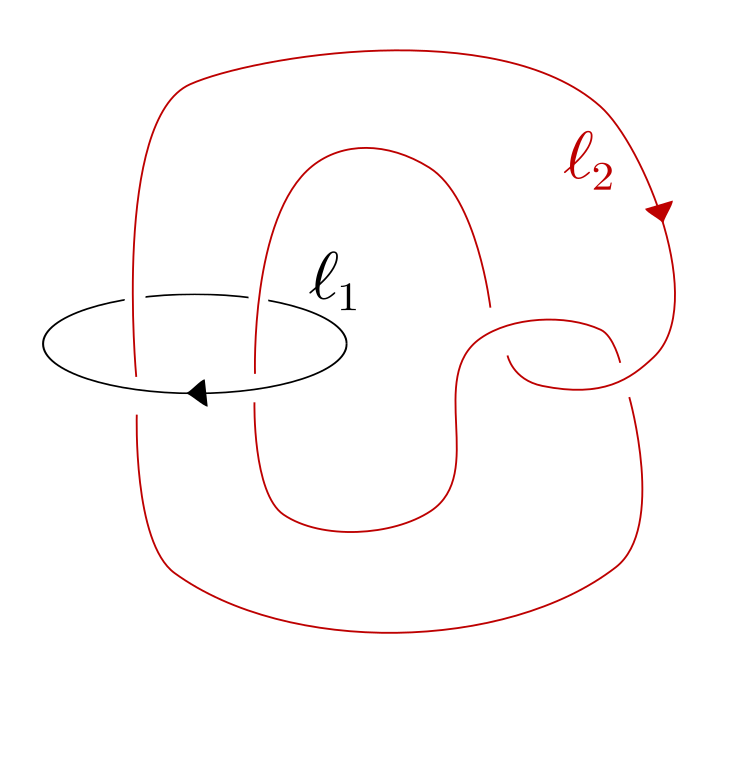}\includegraphics[width=0.45\linewidth]{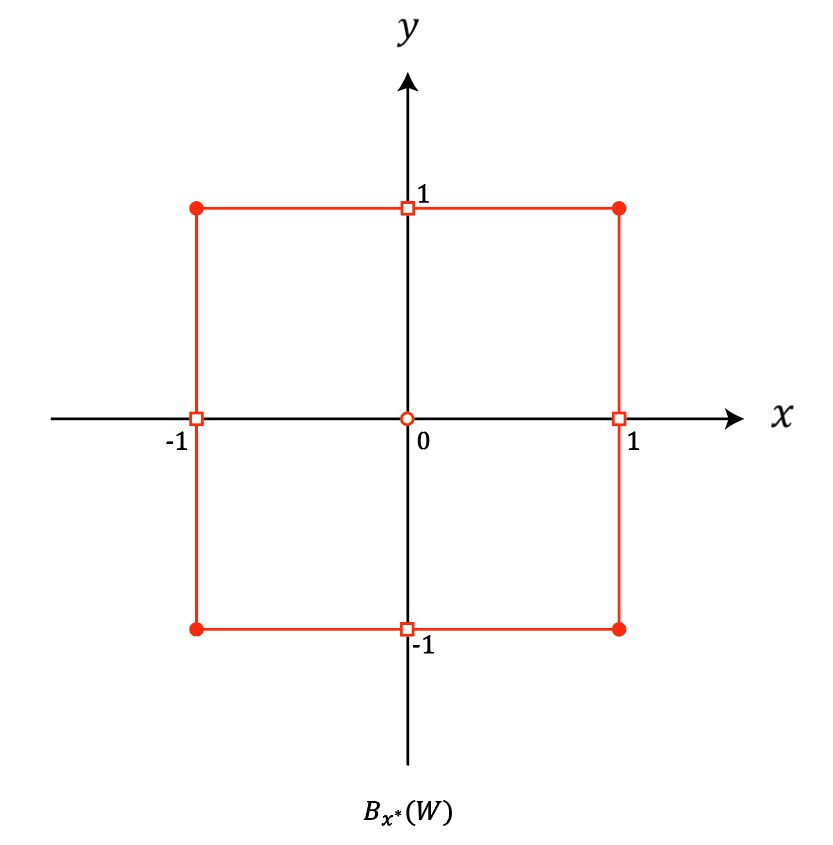}
    \caption{Left: the Whitehead link $L$. Right: the nine integral points with dual Thurston norm at most one in $H^2(M_L,\partial M_L;\R)$. By Gabai, the four vertices are relative Euler classes of some taut foliations whose boundary foliation has no Reeb components. The origin is not the relative Euler class of any taut foliation by Fan--Lai--Yu. Finally, the remaining hollow square points are relative Euler classes of taut foliations with Reeb components at the boundary. \textit{Image credit: Fan--Lai--Yu \cite[Figure 2] {fan2025noteeulerclass0}}.}
    \label{fig: TBWHL}
\end{figure}

Let us now restrict our attention to the exterior $M_L$ of the Whitehead link $L\subset S^3$. Fix an orientation for the components $\ell_1$ and $\ell_2$ of $L$, as on the left of Figure \ref{fig: TBWHL}, then let $\mu_1$ and $\mu_2$ be their respective oriented meridians. We can identify $H^2(M_L,\partial M_L;\R)\cong H_1(M_L;\R)\cong  (\R\mu_1)\oplus (\R\mu_2)$. By \cite[Example 1]{Thurston1986ANF}, the unit ball of the dual Thurston norm on $H^2(M_L,\partial M_L;\R)$ is the square $\{a_1\mu_1+a_2\mu_2\,|\; |a_i|\le 1 \text{ for } i=1,2\}$. There are nine integral points in this square, namely the points $a_1\mu_1+a_2\mu_2$ with $a_1,a_2\in\{0,\pm1\}$. 

By Gabai's Theorem \ref{thm: gabai vertex}, the vertices $\pm (\mu_1 +\mu_2)$ and $\pm (\mu_1-\mu_2)$ are relative Euler classes of taut foliations without boundary Reeb components. 

The remaining five points do not satisfy the parity condition. Fan, Lai and Yu showed that the null class is not realisable as the relative Euler class of any taut foliation of $M_L$ \cite{fan2025noteeulerclass0}. Indeed, if such a foliation existed, it would have a nonzero number of Reeb annuli on each boundary component. A classical theorem of Gabai would then imply that many Dehn-fillings of $M_L$ admit taut foliations. More precisely, this would hold for every manifold $N$ obtained after Dehn-filling the boundary components of $M_L$ along slopes different from the slope of the Reeb components. Finally, Fan, Lai and Yu exhibit ``enough" Dehn-fillings of $M_L$ to find a contradiction. A similar approach to that of Fan--Lai--Yu was also used by Schmalian to prove that the Whitehead link exterior admits no Anosov flows \cite{https://doi.org/10.1112/blms.70049}.

\begin{figure}
    \centering
    \includegraphics[width=0.47\linewidth]{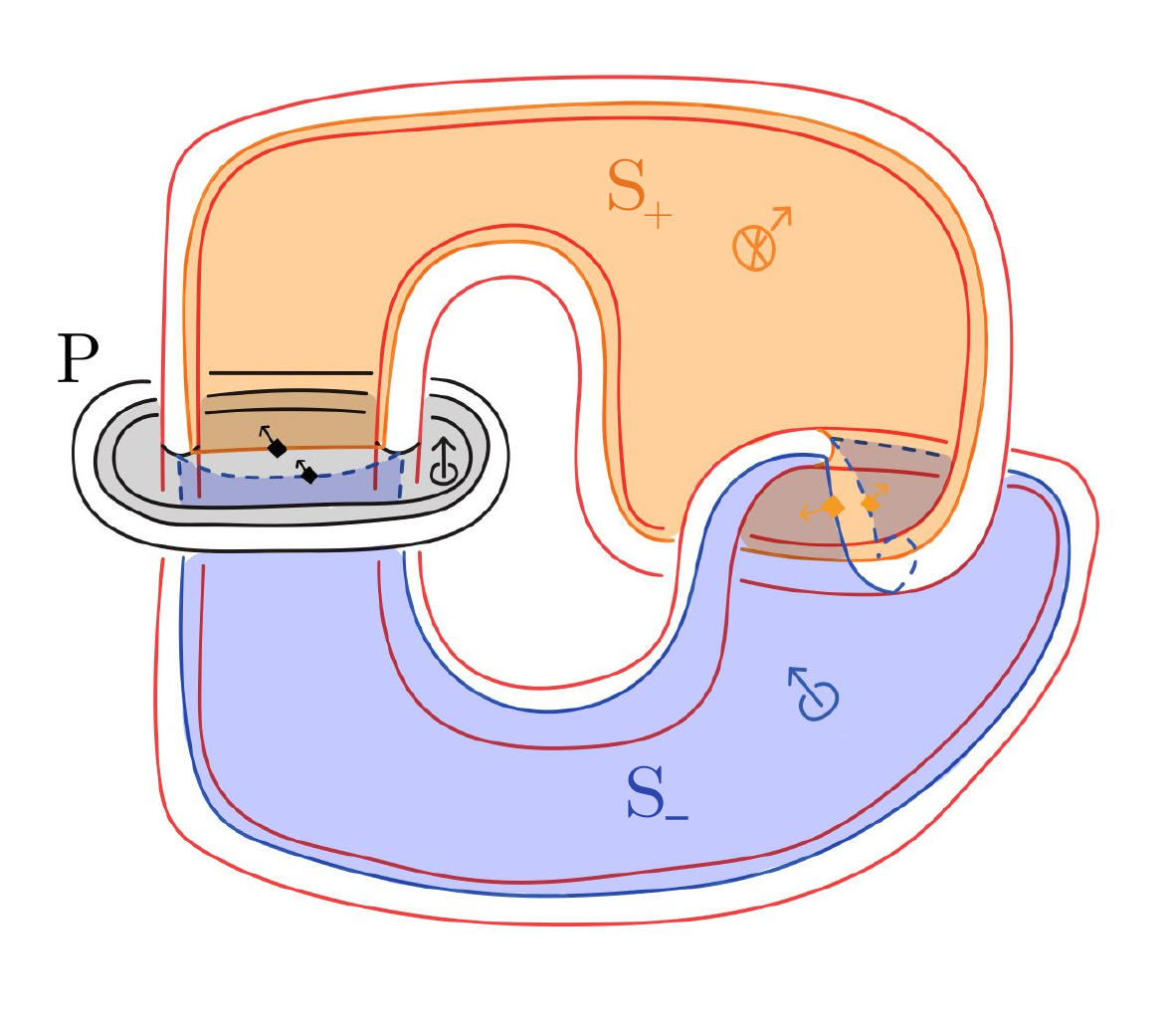}\includegraphics[width=0.5\linewidth]{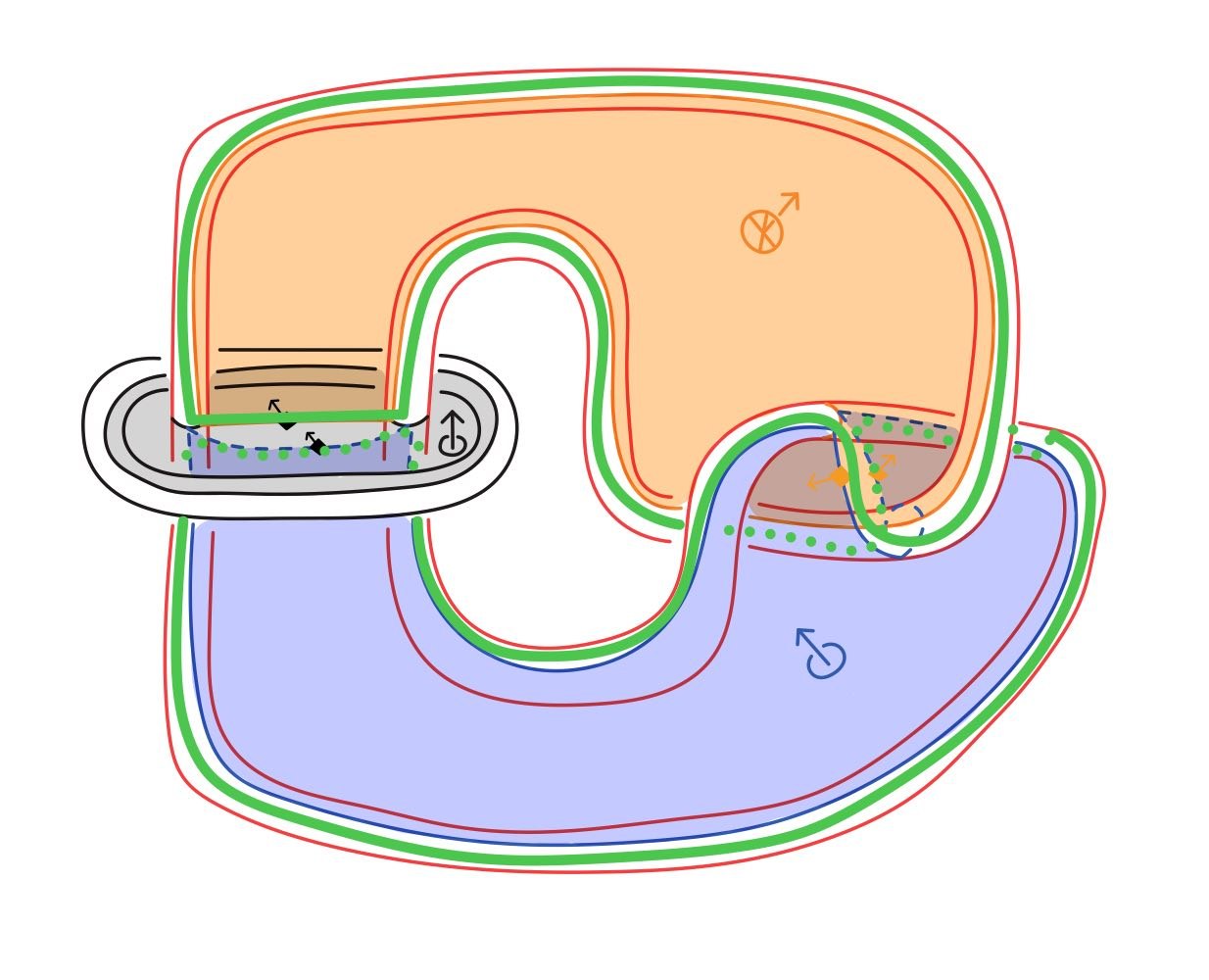}
    \caption{Left: the branched surface $\mc B$. Big arrows indicate the coorientation of the respective sectors. The small arrows indicate the branching directions: these point towards different sectors of $P$ on $P\cap\brl(\mc B)$ and towards the same sector $s_+^2$ of $S_+$ on $S_+\cap\brl(\mc B)$. Right: the green curves outline the components of $\gamma(\mc B)$ in the exterior solid torus $M(\mc B)$.}
    \label{fig: whitehead link branched surface}
\end{figure}

We can now complete the picture above by showing:

\bt\label{thm: whitehead} The relative Euler classes of taut foliations in the Whitehead link exterior are exactly the integral points of unit dual Thurston norm in $H^2(M_L,\partial M_L;\R)$, as on the right of Figure \ref{fig: TBWHL}.
\et
\bp After Gabai and Fan--Lai--Yu's results, we are left to show that the four classes $\pm \mu_1$ and $\pm \mu_2$ are relative Euler classes of taut foliations of the Whitehead link exterior $M_L$. By the symmetry between $\ell_1$ and $\ell_2$, it is enough to construct a taut foliation with relative Euler class $\mu_1$. We will first define a cooriented branched surface $\mc B$ in $M_L$, as depicted on the left of Figure \ref{fig: whitehead link branched surface}. Then, we will show that $\mc B$ fully carries a taut foliation. Finally, we will construct the maw dual graph $\Gamma_m(\mc B)$ to verify it represents the class $\mu_1$.

The component $\ell_1$ of $L$ is unknotted, hence bounds a disk $\Delta$ in $S^3$. The disk $\Delta$ can be isotoped to intersect $\ell_2$ twice and with opposite sign. Let $P$ be the pair of pants $\Delta\cap M_L$, oriented so to intersect $\mu_1$ negatively. The boundary of $P$ divides the boundary torus $\partial N(\ell_2)$ into two annuli $A_+$ and $A_-$. For $\epsilon=+,-$, we define an embedded annulus $S_\epsilon\subset M_L$ with $\partial S_\epsilon\subset \partial M_L\cup P$. One component of $\partial S_\epsilon$ is a meridian of $\ell_1$, and the other component is the union of an essential arc in $P$ and an essential arc in $A_\epsilon$. Orient $S_+$ and $S_-$ so that $S_+$ intersects $\mu_2$ positively and $S_-$ intersects $\mu_2$ negatively. 

The branched surface $\mc B$ is obtained from $P\cup S_+\cup S_-$ as follows. Smooth $S_+$ close to $P\cap S_+$ according to the local coorientations to define a branched surface $\mc B_0$ from $P\cup S_+$. The annulus $S_-$ is cut into a disk by $\mc B_0$. Isotope $S_--\mc B_0$ close to $S_-\cap\mc B_0$ to obtain a cooriented branched surface $\mc B$. 

The exterior $M(\mc B)$ of $\mc B$ is a solid torus, and a meridional disk $D$ intersects $\gamma (\mc B)$ twice (Figure \ref{fig: whitehead link branched surface} right). Hence, the exterior of $\mc B$ in $M_L$ is a product solid torus $(A\times [0,1],\partial A\times [0,1])$, where $A$ indicates a standard annulus. Since the branching locus of $\mc B$ has no triple points, $\mc B$ fully carries a lamination. This lamination can be extended to a foliation $\mc F$ of $M_L$. Furthermore, since $M(\mc B)$ is connected, every sector of $\mc B$ intersects a positive closed transversal, hence $\mc F$ is taut.

\begin{figure}
    \centering
    \includegraphics[width=0.6\linewidth]{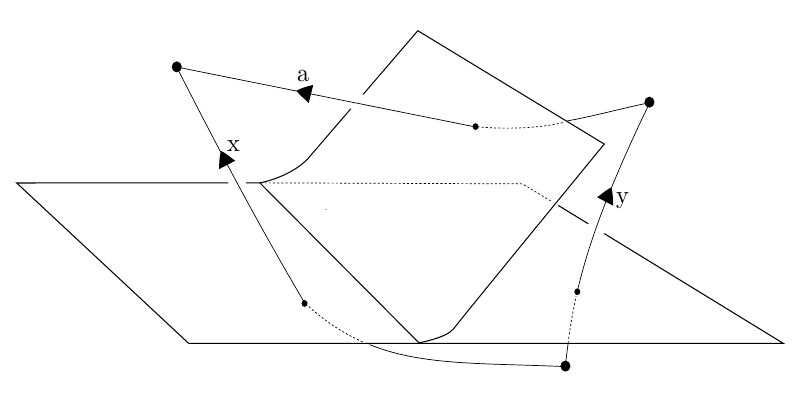}
    \caption{If three sectors of $\mc B'$ share an edge of the branching locus, their dual arcs are linearly dependent. In the case of the picture, equality $y+a=x$ holds in $H_1(M_L)$. Since $a_{\mc B'}(D)$ is null-homologous, if $a=a_{\mc B'}(D)$ then $x$ and $y$ are homologous.}
    \label{fig: null sector}
\end{figure}

Unfortunately, we cannot directly apply the maw dual graph construction to $\mc B$ because its complement is not a union of balls. However, a simplification takes place in the computations, and the relative Euler class $e(\mc F)$ can be read from $\mc B$ anyway. Consider again the meridional disk $D\subset M(\mc B)$ that intersects $\gamma(\mc B)$ twice. Since the foliation $\mc F$ restricts to the product foliation on $M(\mc B)$, the disk $D$ can first be isotoped so that the induced foliation $\mc F\cap D$ is a foliation by intervals $I$ on $D\cong I\times [0,1]$. For any $\varepsilon >0$, one can then isotope $\mc F$ in a neighbourhood of $D$ so that $I\times(\varepsilon,1-\varepsilon)$ is contained in a leaf of $\mc F$, and $I\times[0,\varepsilon]\cup[1-\varepsilon,1]$ is still foliated by intervals. In particular, one can smooth $D$ close to $\partial D$ to obtain a branched surface $\mc B'$ from $\mc B\cup D$ so that $\mc B'$ still fully carries $\mc F$, and its exterior is a product ball. The relative Euler class $e(\mc F)\in H^2(M_L,\partial M_L)$ is now Poincar\'e dual to $[\Gamma_m(\mc B')]\in H_1(M_L)$. Since the complement of $\mc B'$ is connected, each arc $a_{\mc B'}(s)$ dual to a sector $s$ of $\mc B'$ gives rise to a simple closed curve and $$[\Gamma_m(\mc B')]=\sum_{s \text{ sector of }\mc B'} \chi_m(s)[a_{\mc B'}(s)].$$
The anticipated simplification is that $a_{\mc B'}(D)$ vanishes. Indeed, $a_{\mc B'}(D)$ is homotopic to a longitude of $\ell_1$ and, as such, is zero in $H_1(M_L)$. Now, if two sectors $s_1$ and $s_2$ of $\mc B'$ belong to the same sector of $\mc B$, then $a_{\mc B'}(s_1)$ and $a_{\mc B'}(s_2)$ are homologous. This is because the sector $D$ on $\mc B'$ is a source (the branching direction always points out of it). See Figure \ref{fig: null sector}. In particular, \begin{align*}
[\Gamma_m(\mc B')]&=\sum_{s \text{ sector of }\mc B'} \chi_m(s)[a_{\mc B'}(s)]= \sum_{s \text{ sector of } \mc B} \chi_m(s)[a_{\mc B}(s)], \end{align*}
where $a_{\mc B}(s)$ is any simple closed curve in $M_L$ intersecting $\mc B$ exactly once, positively and in $s$. 

The branched surface $\mc B$ has five sectors: two sectors $p_1,p_2$ patch up to $P$, two other sectors $s_+^1,s_+^2$ patch up to $S_+$, and the last one is induced by $S_-$. Choose $s^1_+$ to be the sector of $S_+$ disjoint from $P$, and $p_1$ to be the sector of $P$ disjoint from $\partial N(\ell_1)$. A more attentive inspection allows the reader to fill in Table \ref{table}. Notice that a choice of $a_{\mc B}(s)$ is simple for the sectors of $p_2, s_+^2$ and $S_-$: for each of them, there is an oriented meridian $\mu_i$ intersecting it once and transversely. For the remaining sectors, the relation of Figure \ref{fig: null sector} might be useful.

By adding up the contributions on the right-most column of Table \ref{table}, we obtain $e(\mc F)=[\Gamma_m(\mc B)]=\mu_1$, as desired. \ep

\begin{table}
\centering
\begin{tabular}{c c c c}
\toprule
$s$   &  $\chi_m(s)$ & $[a(s)]\in H_1(M)$ & total contribution \\ \midrule
$p_1$      & $0$        & $-\mu_1-\mu_2$              &    $0$                \\ 
$p_2$      & $-1$        & $-\mu_1$              &    $\mu_1$                              \\
$s_+^1$      & $+1$        & $2\mu_2$              &    $2\mu_2$ \\
$s_+^2$      & $-1$        & $\mu_2$              &    $-\mu_2$     \\
$S_-$       & $+1$        & $-\mu_2$              &    $-\mu_2$                \\\bottomrule
\end{tabular}\caption{Contribution of the sectors of $\mc B$ to $e(\mc F).$ Here, the normal orientation of $\partial M_L$ points out of $M_L$. }\label{table}
\end{table}

\section{Previous work on graphs and taut foliations}\label{sec: previous work}

In the present section, we overview some results by Lackenby \cite{Lackenby2000} and Dunfield \cite{dunfield_floer_2020}. In both works, graphs were used to play the role of Euler classes of taut foliation. We will explain how our construction is related to theirs.

%\vspace{0.1cm}

\subsection{Taut ideal triangulations of $3$-manifolds} 

\begin{figure}
    \centering
    \includegraphics[width=\linewidth]{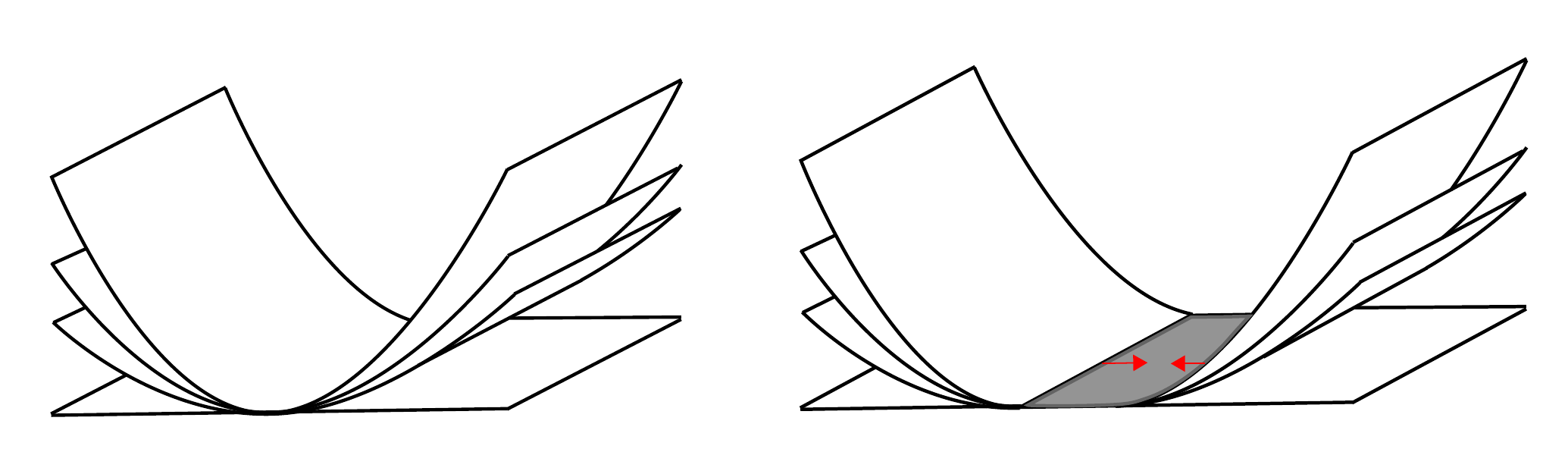}
    \caption{Left: local view of the branched surface carried by a taut triangulation close to an edge. Right: the resulting branched surface after a flattening operation. A new rectangular sector is introduced (dark grey in the picture). The maw vector field is now well-defined on the branched surface, as indicated by the red arrows.}
    \label{fig: flattening}
\end{figure}

In \cite{Lackenby2000}, Lackenby introduced the notion of \emph{taut triangulations} for a $3$-manifold $M$ with nonempty boundary. Roughly speaking, a taut triangulation of $M$ is a triangulation of $M-\partial M$ by ideal $3$-simplices with a special angled structure. Two faces of each $3$-simplex have normal orientation pointing inward, and the other two have normal orientation pointing outward. An edge between two coherently cooriented faces has interior angle $\pi$, otherwise it has zero interior angle. It is understood that the glueing map between two taut $3$-simplices must respect their normal orientations. Moreover, the total angle around every edge of the triangulation should be $2\pi$. See Figure \ref{fig: flattening} left for a local picture of $M$ close to an edge of the triangulation. 

The $2$-skeleton of a taut triangulation yields a branched surface $\mc B$ in $M$. Unfortunately, the branching locus of $\mc B$ is not \emph{generic}, i.e. it is not modelled as in Figure \ref{fig: branching models}. However, as no triple points arise, the branched surface fully carries a lamination $\mc L$. Notice that we can flatten $\mc B$ out close to each edge by introducing a new rectangular sector (Figure \ref{fig: flattening}). After this procedure, the new branched surface still fully carries $\mc L$, and a maw vector field is now well-defined.

The cusped structure of taut simplices makes the exterior of $\mc B$ a union of product balls. In particular, $\mc B$ carries a cooriented foliation that is taut if $M$ is compact (e.g., the graph $G$ defined below can be split into a number of closed transversals intersecting every leaf of the foliation).

For this reason, taut ideal triangulations were employed by Lackenby to give foliation-free proofs of classical results on manifolds admitting taut foliations. For instance, it was shown in \cite{Lackenby2000} that manifolds admitting taut triangulations are irreducible and $\partial$-irreducible. Moreover, a new proof of Gabai's Theorem that the minimal genus of a knot equals its singular genus was presented.

\bigskip

In order to substitute the Euler class of the foliation carried by a taut triangulation with a combinatorial analogue, Lackenby uses the graph $G$ in $M$ dual to the taut triangulation. Every edge of $G$ inherits the normal orientation of the $2$-simplices of the triangulation. Since every vertex of $G$ has exactly two incoming edges and exactly two outgoing ones, $G$ is a cycle and hence represents an element of $H_1(M)$. Through a study of the singular normal surfaces carried by the triangulation, Lackenby shows:

\bprop\cite[Proposition 11]{Lackenby2000}\label{prop: lackenby} Let $M$ be a compact $3$-manifold with a taut triangulation and let $G$ be the graph described above. For every compact oriented surface $F$ without disk or sphere components, and for every continuous map $(F,\partial F)\to (M,\partial M)$, we have $$|\langle[G],[F]\rangle|\le -2\chi(F).$$
\eprop

\begin{figure}
    \centering
    \includegraphics[width=0.9\linewidth]{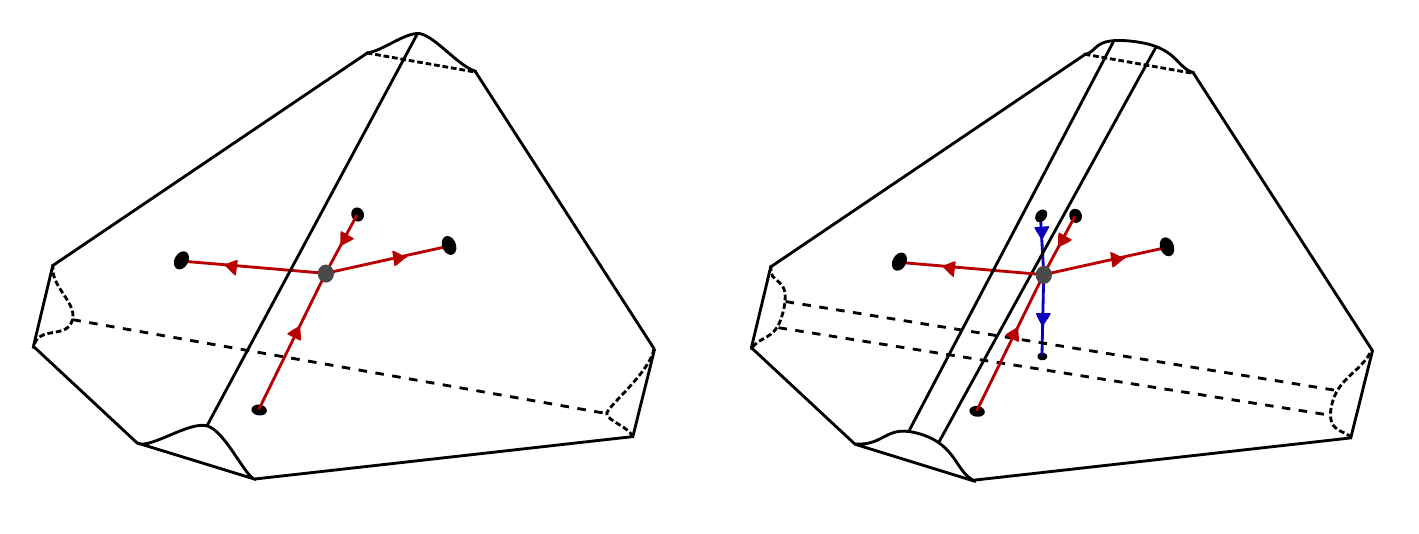}
    \caption{Left: the dual graph $G$ coloured red. Right: the subdivision of the boundary of a $3$-simplex into sectors after the flattening operation. The local arc of the link $\beta$ is represented in blue.}
    \label{fig: graphs}
\end{figure}

As the equality is achieved by a surface carried by the branched surface associated with the triangulation, such a surface minimises the topological complexity among (singular) surface representatives of its homology class. We show that:

\bl\label{lemma: lackenby's graph} Let $M$ be a compact $3$-manifold with a taut triangulation, $\mc F$ a taut foliation carried by it, and $G$ the dual graph described above.
The relative Euler class $e(\mc F)\in H^2(M,\partial M)$ is Poincar\'e dual to $-\frac 12[G]\in H_1(M)$.
\el
\bp First flatten the branched surface induced by the triangulation close to every edge, as in Figure \ref{fig: flattening}. Call $\mc B$ the resulting branched surface. Endow $\partial M$ with the outward-pointing normal orientation. Since $\mc B$ has a well-defined maw vector field, the proofs of Lemma \ref{lemma: graph is cycle} and of Theorem \ref{thm: maw dual graph} can be carried out verbatim. The associated maw dual graph $\Gamma_+$ represents the Poincar\'e dual of $e(\mc F)$. The sectors of $\mc B$ are of two types: hexagons arising as truncated $2$-simplices or rectangles arising from the flattening procedure. The branching direction points outside at every edge of the hexagons and points inside at every edge of the rectangle. As $\partial M$ is normally oriented outwardly, hexagonal sectors have $\chi_m=+1$ whilst rectangular sectors have $\chi_m=-1$. In particular, $$[\Gamma_+]=[G]+[\beta],$$ where $\beta\subset M$ is a link intersecting every rectangular sector once with negative sign, see Figure \ref{fig: graphs}.

Now consider the inward-pointing normal orientation on $\partial M$. This time, hexagonal sectors of $\mc B$ have $6$ double corners, hence $\chi_m=-2$. Rectangular sectors are now sink sectors, and $\chi_m=+1$. The maw dual graph $\Gamma_-$ associated with the new orientation of the boundary satisfies $$[\Gamma_-]=-2[G]-[\beta].$$

Since the two boundary sections are homotopic, they induce the same relative Euler class $e(\mc F)$. Therefore $[\Gamma_+]=[\Gamma_-]$ implies $[\beta]=-\frac 32[G]$ and $$[\Gamma_+]=[G]+[\beta]=[G]-\frac 32 [G]=-\frac 12 [G].$$
\ep

Proposition \ref{prop: lackenby} then gives a new foliation-free proof of Thurston's Inequality for singular surfaces with respect to taut foliations carried by taut ideal triangulations:

\bc\label{Lackenby's Inequality} Let $M$ be a compact manifold admitting a taut triangulation, and let $\mc F$ be a carried foliation. Then for every compact oriented surface $S$ without disk or sphere components, and for every continuous map $(S,\partial S)\to (M,\partial M)$, we have \begin{equation}\label{Thurston's Inequality}
    |\langle e(\mc F),[S]\rangle|\le |\chi(S)|.
\end{equation}
\ec

Thurston first proved expression (\ref{Thurston's Inequality}) for taut foliations and properly embedded surfaces \cite{Thurston1986ANF}. The main idea was to isotope the surface so that only saddle tangencies occurred with the foliation. Gabai extended this proof to singular surfaces, first for finite-depth foliations \cite{gabai_foliations_1983}, then to the general case \cite{Gabai2000}. As Gabai explains in the introduction of \cite{Gabai2000}, for $C^2$ taut foliations, a minimal surface proof of the most general result was already known by putting together works of Sullivan \cite{Sullivan1979}, Schoen and Yau \cite{SchoenYau1979} and Hass \cite{Hass1986}.

\vspace{0.1cm}

\subsection{Foliar orientations} 

In \cite[Section 7]{dunfield_floer_2020}, Dunfield used a special kind of triangulation to construct taut foliations on closed $3$-manifolds. Given a triangulation $\mc T$ of a closed $3$-manifold $M$, an \emph{edge orientation} is a choice of orientation on the edges of $\mc T$. An edge orientation is \emph{acyclic} if no $2$-simplex of $\mc T$ has a directed cycle as boundary orientation. In particular, every triangle has two \emph{short edges} and one \emph{long edge}, as in \cite[Figure 6]{dunfield_floer_2020}. The existence of an acyclic edge orientation guarantees we can define a cooriented branched surface $\mc B$ in $M$ dual to the triangulation, as in \cite[Figure 7]{dunfield_floer_2020}. An acyclic edge orientation is \emph{foliar} if the associated branched surface is laminar in the sense of Li \cite{li_laminar_2002} and its exterior consists of product balls. In fact, the property of being foliar can be completely described just in terms of the edge orientation, i.e., without referring to the branched surface \cite{dunfield_floer_2020}. It follows that if $\mc T$ admits a foliar orientation, then $M$ has a taut foliation.

If the manifold $M$ has a taut foliation with vanishing Euler class, then its fundamental group is left-orderable, and this is of significant interest for the L-space Conjecture \cite{dunfield_floer_2020}. Dunfield then shows how to compute the Euler class of a taut foliation induced by a foliar orientation. This is now explained. 

Suppose that the triangulation $\mc T$ on $M$ is endowed with a foliar orientation. An edge $\epsilon$ of $\mc T$ is \emph{mixed} for the $3$-simplex $\Delta$ if $\epsilon$ appears on the $1$-skeleton of $\Delta$ as the long edge of a face and a short edge of another face.

Since the sectors of $\mc B$ generate the simplicial $2$-chains of $M$, we can define a simplicial $2$-cochain $\phi$ by assigning an integer to each of them. For every sector $D$ of $\mc B$, we put $$\phi(D):=1-\frac 12 \text{mixed}(\epsilon),$$ where $\epsilon$ is the edge of $\mc T$ dual to $D$, and $\text{mixed}(\epsilon)$ counts the number of times $\epsilon$ appears as a mixed edge in every tetrahedron of $\tau$. After this definition, Dunfield shows: 

\bt\cite[Theorem 9.2]{dunfield_floer_2020} Suppose $\mc F$ is the foliation of $M$ associated to a foliar orientation of a $1$-vertex triangulation $\mc T$. Then the cochain $\phi$ is a cocycle representing the Euler class $e(\mc F)\in H^2(M)$.
\et

The original proof of this result consists in defining a vector field on $\mc B$ to compute $e(\mc F)$. Alternatively, our formula for the Euler class can be used too, even after removing the assumption that the triangulation has exactly one vertex. Indeed, given a sector $D$, the triple point on $\partial D$ inside a tetrahedron $\Delta$ gives rise to a double corner for $D$ if and only if the dual edge $\epsilon$ of $D$ is mixed for $\Delta$, see Figure 7 of \cite{dunfield_floer_2020}. Since $D$ is a disk, $1-\frac 12 \text{mixed}(\epsilon)$ is exactly $\chi_m(D)$.

\section{A Transverse Surface Theorem for branched surfaces}\label{sec: carried vs fully marked}
Let $\mc B$ be a branched surface in a compact oriented $3$-manifold $M$. Suppose that the complement of $\mc B$ in $M$ is a union of balls. It is natural to ask which topological information is retained by $\mc B$. For example, Landry \cite{LANDRY} and Landry--Minsky--Taylor \cite{LandryMinskyTaylor2024, LandryMinskyTaylor+2026+203+257} studied the intertwining between surfaces in $M$ and $\mc B$, when $\mc B$ is the $2$-skeleton of a veering triangulation on $M$. \emph{Veering triangulations} and \emph{strict} veering triangulations are taut ideal triangulations satisfying some special combinatorial properties; for their definition, the reader can consult \cite{LandryMinskyTaylor+2026+203+257}. We are interested in understanding which properties of veering triangulations generalise to our setting.

First, we observe that the directed dual graph $\Gamma(\mc B)$ -- i.e., the unweighted maw dual graph -- determines the first homology of $M$. The following lemma is a generalisation of \cite[Lemma 5.8]{LandryMinskyTaylor2024}. We say that an oriented graph $\Gamma$ is \emph{strongly connected} if for any couple of vertices $u,v$ there is an oriented path on $\Gamma$ walking from $u$ to $v$.

\bl Let $M$ be a compact oriented $3$-manifold and $\mc B\subset M$ a cooriented branched surface whose complement is a union of balls. The inclusion induced map $\iota_*:\pi_1(\Gamma(\mc B))\to \pi_1(M)$ is surjective. Therefore, the (non-directed) cycles on $\Gamma(\mc B)$ generate $H_1(M)$.

Moreover, if $\mc B$ carries a taut foliation, then every connected component of $\Gamma(\mc B)$ is strongly connected and the directed cycles generate $H_1(M)$.
\el
\bp Let $\gamma\subset M$ be a closed loop. We want to homotope $\gamma$ to lie on $\Gamma(\mc B)$. First, perform a homotopy so that $\gamma$ is transverse to $\mc B$ and does not intersect the branching locus. After another homotopy, we can suppose that every intersection of $\gamma\cap \mc B$ happens at the base point of a sector of $\mc B$. Now, let $D$ be a complementary ball of $\mc B$ and consider a component $\alpha$ of $\gamma\cap D$. If $\alpha$ is a closed component, then it is null-homotopic because it is contained in a ball. Otherwise, $\alpha$ has endpoints $x$ and $y$ and there is a shortest path $\delta\subset \Gamma(\mc B)$ with the same endpoints. The curve $\alpha\cup\delta$ is a closed loop in a ball, hence null-homotopic. So, we can homotope $\alpha$ to $\delta$ rel. $x$ and $y$. After repeating this reasoning for every component $\alpha$ and for every complementary ball $D$, we obtain what we wanted.

\bigskip

Suppose now that $\mc B$ carries a taut foliation. Let $u$ be a vertex of $\Gamma(\mc B)$ and let $S(u)$ be the closure of the union of the complementary balls of $\mc B$ whose base points are reachable from $u$ through a directed path on $\Gamma(\mc B)$. The boundary of $S(u)$ is a union of portions of $\partial M$ and of sectors of $\mc B$. No such sector could have normal orientation pointing out of $S(u)$, because otherwise there would be a complementary ball not in $S(u)$ although reachable from $u$ through a directed path. In particular, the normal orientation of $\partial S(u)\cap \mc B$ points always inside $S(u)$. Since $\mc B$ carries a taut foliation, $\partial S(u)\cap \mc B=\emptyset$ and $S(u)$ must be a connected component of $M$. This shows that each component of $\Gamma(\mc B)$ is strongly connected. The fact that directed cycles generate follows as in the proof of \cite[Lemma 5.8]{LandryMinskyTaylor2024}.
\ep

We now focus on the Combinatorial Transverse Surface Theorem (TST).

\bt[Combinatorial TST]\cite[Theorems 5.4 and 5.5]{LANDRY} Let $M$ be a compact oriented $3$-manifold and $\tau$ a veering triangulation of it. An integral homology class $\alpha\in H_2(M,\partial M)$ pairs nonnegatively with every directed cycle in $\Gamma(\tau^{(2)})$ if and only if a surface representative of $\alpha$ is carried by $\tau^{(2)}$. Moreover, such a surface has no null-homologous components and realises the Thurston norm of $\alpha$.
\et

We observe that the proof of this previous result easily extends to the context of cooriented branched surfaces whose complement consists of balls.

\bt[Combinatorial TST for branched surfaces]\label{thm: TST} Let $M$ be a compact oriented $3$-manifold and $\mc B\subset M$ a cooriented branched surface with generic branched locus and complement a union of balls. An integral homology class $\alpha\in H_2(M,\partial M)$ pairs nonnegatively with every directed cycle in $\Gamma(\mc B)$ if and only if a surface representative $S$ of $\alpha$ is carried by $\mc B$. Moreover, $S$ has no null-homologous components and, if $\mc B$ carries a taut foliation $\mc F$, $S$  realises the Thurston norm of $\alpha$.
\et
\bp Since $\Gamma(\mc B)$ intersects positively each sector of $\mc B$, any surface $S$ carried by $\mc B$ must pair nonnegatively with directed cycles in $\Gamma(\mc B)$. In particular, no component of $S$ can be nullhomologous. Furthermore, as $S$ can be written as a glueing of copies of sectors of $\mc B$, it satisfies $$\langle [\Gamma_m(\mc B)],[S]\rangle=\chi(S).$$ Now, if $\mc B$ carries a taut foliation, Thurston's Inequality guarantees that $S$ is norm-minimising.

\begin{figure}
    \centering
    \includegraphics[width=0.7\linewidth]{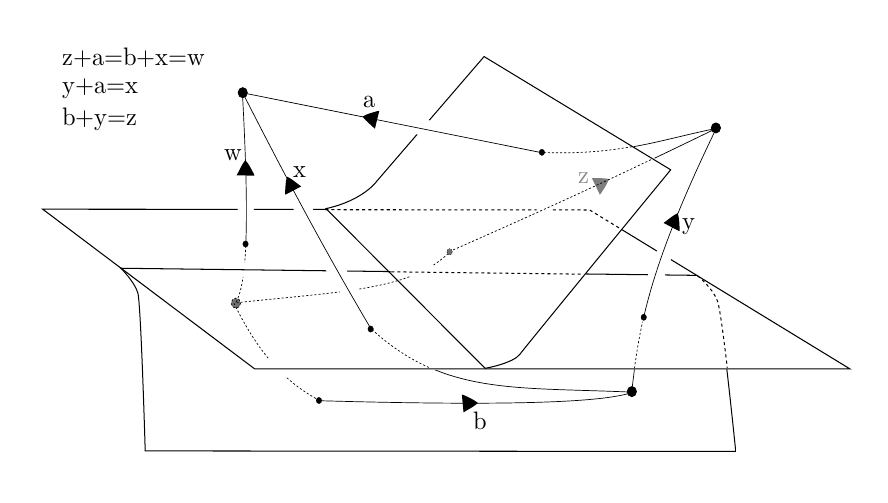}
    \caption{Local picture of the branched surface $\mc B$ close to the branching locus. The labels $a,b,x,y,z,w$ are the integers assigned by $m$ to the respective arc and by $\sigma$ to the respective dual sector. Since the unoriented arcs labelled $x,y$ and $a$ patch to a cycle $\gamma\in \Gamma(\mc B)$ that bounds a disk in $M$, we have $y+a-x=\langle i^*\eta,\gamma\rangle=\langle\eta,i_*\gamma\rangle=0$. Similarly, the branching equations on the top-left of the picture hold. For every sector $s\in \mc B$, take $\sigma(s)$ parallel copies of it. The branching equations guarantee that the copies of all sectors finally patch together to a properly embedded surface $S$ carried by $\mc B$.}
    \label{fig: branching equations}
\end{figure}

Conversely, suppose that a certain homology class $\alpha\in H_2(M,\partial M)$ pairs nonnegatively with all directed cycles in $\Gamma(\mc B)$. Let $\eta\in H^1(M)$ be the Poincar\'e dual of $\alpha$ and $i:\Gamma(\mc B)\to M$ the inclusion map. The induced class $i^*\eta\in H^1(\Gamma(\mc B))$ pairs nonnegatively with all directed cycles in $\Gamma(\mc B)$. By \cite[Lemma 5.10]{LandryMinskyTaylor2024}, there is a nonnegative $1$-cocycle $m:\{\text{edges of }\Gamma(\mc B)\}\to\R$ representing $i^*\eta$. Via the identification between sectors of $\mc B$ and edges of $\Gamma(\mc B)$, we obtain a nonnegative map $\sigma: \{\text{sectors of } \mc B\}\to \R$. Since $m$ represents $i^*\eta$, $\sigma$ is a combinatorial measure on the branched surface $\mc B$. Namely, at any point of the branching locus, $\sigma$ satisfies the branching equations (see Figure \ref{fig: branching equations}). In particular, there is a surface $S$ carried by $\mc B$ and such that the algebraic intersections of $\alpha$ and $[S]$ with any cycle of $\Gamma(\mc B)$ coincide. Since the cycles of $\Gamma(\mc B)$ generate $H_1(M)$, the surface $S$ represents $\alpha\in H_2(M,\partial M)$. 
\ep

Landry, Minsky and Taylor showed that a stronger TST holds when $\mc B$ is the $2$-skeleton of a \emph{strict} veering branched surface \cite[Theorem 2.13]{LandryMinskyTaylor+2026+203+257}. 

\bt\cite[Combinatorial Strong TST]{LandryMinskyTaylor+2026+203+257}\label{thm: strong TST} Let $\tau$ be a strict veering triangulation of $M$. Let $S$ be a norm-minimising surface in $M$ without null-homologous components. The following are equivalent:
\begin{itemize}
    \item[(a)] $[S]$ pairs nonnegatively with every directed loop in $\Gamma(\tau^{(2)})$;
    \item[(b)] $x([S])=-e_\tau([S])$;
    \item[(c)] $S$ is carried by $\tau$ up to an isotopy, which can be chosen to fix $\partial S$. 
\end{itemize}
\et

Here, the \emph{Euler class} $e_\tau$ of the veering triangulation $\tau$ is defined as in \cite[Section 2.11]{LandryMinskyTaylor+2026+203+257}. Since we are considering non-relative veering triangulations, the definition of $e_\tau$ is $-\frac 12$ of the one for taut ideal triangulations as defined by Lackenby. By Lemma \ref{lemma: lackenby's graph}, $e_\tau$ is the Euler class of any foliation carried by the $2$-skeleton of the triangulation, hence represented by $\Gamma_m(\tau^{(2)})$.

Unfortunately, there is no hope of straightforwardly extending the Combinatorial Strong TST to branched surfaces whose complement is a union of balls: let us construct a counterexample. Let $\Sigma$ be an oriented connected surface of genus one with one boundary component, and let $M:=\Sigma\times S^1$ be endowed with the product foliation $\mc F$. Let $\alpha,\beta\subset \Sigma$ be two disjoint properly embedded arcs such that $\Sigma-(\alpha\cup\beta)$ is homeomorphic to a disk. Fixed a fiber $\Sigma_0:=\Sigma\times \{\star\}\subset M$, the annulus $\alpha\times S^1$ is cut into a disk $D_\alpha$ by $\Sigma_0$. Similarly, $\beta\times S^1$ traces a disk $D_\beta\subset M-\Sigma_0$. Choose normal orientations for $D_\alpha$ and $D_\beta$. The foliation $\mc F$ is carried by the branched surface $\mc B$ obtained by smoothing $\Sigma_0\cup D_\alpha \cup D_\beta$ according to the normal orientations. We can require $\mc B$ to carry the oriented cut and paste $S$ of $\Sigma_0$ and $\alpha\times S^1$. In this way, equality $\langle e(\mc F),[S]\rangle=\chi(S)$ holds. The counterexample is given by considering the branched surface $\mc B'$ obtained by smoothing $\Sigma_0\cup\overline{D_\alpha}\cup D_\beta$ according to coorientations. The branched surface $\mc B'$ still carries $\mc F$, but now $S$ intersects negatively the arc of $\Gamma(\mc B')$ dual to $\overline{D_\alpha}$.

\bibliographystyle{alpha}
\bibliography{bibliography}

\end{document}